\documentclass[11pt,a4paper]{article}
\usepackage{pifont}
\usepackage{CJK,epsfig,amsfonts,amsgen,amsmath,amstext,amsbsy,amsopn,amsthm
}
\usepackage{mathrsfs}
\usepackage{float}
\usepackage{epsfig,latexsym,amsfonts,amsmath,amssymb}
\usepackage{color}
\usepackage{calc}
\usepackage{enumerate}

\usepackage{comment}
\usepackage{geometry}
\usepackage{calc}
\usepackage{multirow}

\usepackage{tikz}
\usetikzlibrary{shapes.geometric, arrows, positioning}
\usetikzlibrary{decorations.markings}
\tikzstyle{vertex}=[circle, draw, inner sep=0pt, minimum size=6pt]
\tikzstyle{node} = [rectangle, minimum width=1.5cm, minimum height=0.7cm, text centered, draw=black, fill=blue!20, font=\small]
\tikzstyle{arrow} = [thick,->,>=stealth]

\newtheorem{theorem}{Theorem}

\newtheorem{lemma}[theorem]{Lemma}
\newtheorem{conjecture}[theorem]{Conjecture}

\newtheorem{case}[]{Case}
\newtheorem{claim}[]{Claim}
\newtheorem{subcase}{Subcase}[case]

\newtheoremstyle{myremark}%
  {}{}%
  {\upshape}{}
  {\bfseries}{.}%
  { }%
  {}
\theoremstyle{myremark}
\newtheorem{remark}[theorem]{Remark}

\newtheoremstyle{myexample}%
  {}{}%
  {\upshape}{}
  {\bfseries}{.}%
  { }%
  {}
\theoremstyle{myexample}
\newtheorem{example}[theorem]{Example}

\usepackage{indentfirst} 
\newcounter{Hcase}

\theoremstyle{plain}

\date{\today}
\begin{document}

\title{Beyond Frankl's conjecture: Maximum sum of sizes of cross-intersecting families with restricted conditions}

\author{\small Zhen Jia${\thanks{Research partially supported by the National Key R\&D Program of China (No. 2025YFA1017700). E-mail address: 12231271@mail.sustech.edu.cn}}^{~1}$,~Qing Xiang${\thanks{Research partially supported by the National Key R\&D Program of China (No. 2025YFA1017700) and the National Natural Science Foundation of China Grant (No. 12131011). E-mail address: xiangq@sustech.edu.cn}}^{~1}$,  Jimeng Xiao${\thanks{Research partially supported by the National Natural Science Foundation of China Grant (No. 12501460). E-mail address: xiaojm@sustech.edu.cn}}^{~1}$ and Huajun Zhang${\thanks{Research partially supported by the National Natural Science Foundation of China Grant (No. 12371332). E-mail address: huajunzhang@usx.edu.cn}}^{~2}$
 \\[2mm]
\small ${}^1$  Department of Mathematics\\[-0.8ex]
\small Southern University of Science and Technology\\[-0.8ex]
\small Shenzhen, Guangdong 518055, China\\
\small ${}^2$Department of Mathematics\\[-0.8ex]
\small Shaoxing University \\[-0.8ex]
\small Shaoxing, Zhejiang 312000, China}

\baselineskip 20pt

\date{}
\maketitle
 \vspace{4mm}

\begin{abstract}
Let $\mathcal{F}$ and $\mathcal{G}$ be two families of subsets of $[n]$. 
We say that $\mathcal{F}$ and $\mathcal{G}$ are \emph{cross-intersecting} if $F \cap G \neq \emptyset$ 
for all $F \in \mathcal{F}$ and $G \in \mathcal{G}$. 
Let $n, k_1, k_2, m, r, t$ be positive integers satisfying $n  > k_1 + k_2$. 
Suppose that $\mathcal{F} \subseteq \binom{[n]}{k_1}$ and $\mathcal{G} \subseteq \binom{[n]}{k_2}$ 
    are two cross-intersecting families such that 
    $\left\{F \in \binom{[n]}{k_1} : |F \cap [m]| \geq r\right\} \subseteq \mathcal{F}$ 
    and $\mathcal{F}$ is $t$-intersecting. In this paper, we determine the maximum sum $|\mathcal{F}| + |\mathcal{G}|$ under these assumptions. 
By taking $k_1=k_2+t-1$, $m = k_2+t+s-1$ and $r=k_2+t-1$ for $s\ge 0$, our result recovers a conjecture of Frankl \cite{frankl2016new}, which was recently confirmed by Wu, Feng, and Li \cite{wu2025proof}. Further, we establish a complete characterization of all extremal families.
Our approach extends the Ahlswede--Khachatrian generating set method to 
handle the simultaneous constraints of cross-intersection and $t$-intersection, 
subject to additional structural conditions.



\vskip 0.1in \noindent%
 \\[7pt]
\end{abstract}

\maketitle
\section{Introduction}

In this paper, we consider cross-intersecting families. First, we fix some 
notation. Throughout the paper, lowercase letters are used to denote integers, 
uppercase letters for sets of integers, and calligraphic uppercase letters 
for collections of sets. The set of integers from $\ell$ to $n$ inclusive 
is denoted by $[\ell,n]$. In particular, we use $[n]$ as an abbreviation 
for the set $\{1, 2, \ldots, n\}$. A $k$-set (or $k$-subset) is a set of 
cardinality $k$, and the collection of all $k$-subsets of $[n]$ is denoted 
by $\binom{[n]}{k}$.

Let $n$ and $k$ be integers with $1 \leq k \leq n$. A family 
$\mathcal{F} \subseteq \binom{[n]}{k}$ is said to be \emph{$t$-intersecting} 
if $|A \cap B| \geq t$ for all $A, B \in \mathcal{F}$. 
Cross-intersecting families are a natural variant of intersecting families. 
Two families $\mathcal{F}_1 \subseteq \binom{[n]}{k_1}$ and 
$\mathcal{F}_2 \subseteq \binom{[n]}{k_2}$ are said to be 
\emph{cross-intersecting} if $|A \cap B| \geq 1$ for all 
$A \in \mathcal{F}_1$ and $B \in \mathcal{F}_2$.

The celebrated Erd\H{o}s-Ko-Rado theorem \cite{erdos1961intersection} determines the maximum size of intersecting families, and also characterizes the extremal families.
\begin{theorem} {\em \cite{erdos1961intersection}}
    Let $n\geq 2k$. If $\mathcal{F}\subseteq\binom{[n]}{k}$ is an intersecting family, then $|\mathcal{F}|\leq\binom{n-1}{k-1}$. Moreover, if $n>2k$, equality holds if and only if $\mathcal{F}=\{A\in\binom{[n]}{k}:a\in A\}$ for some $a\in [n]$.
\end{theorem}

The Erd\H{o}s--Ko--Rado theorem is widely regarded as a cornerstone result 
of extremal set theory and has many interesting applications and generalizations, see \cite{DF2,ellis2022intersection,frankl2016invitation,frankl1978,hilton1977intersection,hilton1967some,k,k68,k63,W84} 
and the references therein. In particular, $t$-intersecting and 
cross-intersecting problems are two important variations.

One could also consider a further variation by combining the above two 
conditions—that is, requiring two families to be cross-intersecting while 
one of them is additionally $t$-intersecting. In \cite{frankl2016new}, 
Frankl studied such families under the constraint that $\mathcal{F}$ is 
$(t+1)$-intersecting and $|\mathcal{F}| \geq 1$. He established an upper bound 
on sum of their sizes of these families and proposed a more general problem concerning 
the structure when $\binom{[k+t+s]}{k+t} \subseteq \mathcal{F}$.

\begin{theorem}[\cite{frankl2016new}]
    Let $n,k,t$ be non-negative integers, $n\geq 2k+t$. Suppose that $\mathcal{F}\subseteq \binom{[n]}{k+t}$ and $\mathcal{G}\subseteq\binom{[n]}{k}$ are cross-intersecting. 
    \begin{itemize}
        \item[(i)] If $\mathcal{F}$ is $t$-intersecting, then $|\mathcal{F}|+|\mathcal{G}|\leq \binom{n}{k}$. 
        \item[(ii)] If $\mathcal{F}$ is $(t+1)$-intersecting and $|\mathcal{F}|\geq1$, then $|\mathcal{F}|+|\mathcal{G}|\leq 1+\binom{n}{k}-\binom{n-(k+t)}{k}$.
    \end{itemize}
\end{theorem}
\begin{conjecture}[\cite{frankl2016new}]\label{frankl's-conj}
     Let $t\geq1,s\geq0$, $k\geq s+1$ and $n\geq 2k+t-1$ be integers. Suppose that $\mathcal{F}\subseteq \binom{[n]}{k+t-1}$ and $\mathcal{G}\subseteq\binom{[n]}{k}$ are cross-intersecting. If $\mathcal{F}$ is $t$-intersecting and $\binom{[k+t-1+s]}{k+t-1}\subseteq\mathcal{F}$, then
     \begin{equation*}
         |\mathcal{F}|+|\mathcal{G}|\leq \binom{k+t-1+s}{k+t-1}+\binom{n}{k}-\sum_{i=0}^{s}\binom{k+t-1+s}{i}\binom{n-k-t+1-s}{k-i}.
     \end{equation*}
\end{conjecture}
The upper bound on $|\mathcal{F}|+|\mathcal{G}|$ in Conjecture \ref{frankl's-conj} is tight. 
Setting $\mathcal{F}=\binom{[k+t-1+s]}{k+t-1}$ and 
$\mathcal{G}=\left\{G\in\binom{[n]}{k} : |G\cap [k+t-1+s]|\geq s+1\right\}$, 
we see that $\mathcal{F}$ and $\mathcal{G}$ are cross-intersecting, $\mathcal{F}$ is $t$-intersecting, and $|\mathcal{F}|+|\mathcal{G}|$ attains 
this maximum. The recent work of Wu, Feng, and Li \cite{wu2025proof} 
verified this conjecture, successfully establishing the upper bound by 
utilizing Frankl's original induction approach in Conjecture \ref{frankl's-conj}.
 
In the above statement of Frankl's conjecture, the restriction 
$\binom{[k+t+s]}{k+t} \subseteq \mathcal{F}$ can be equivalently 
rephrased as $\left\{F \in \binom{[n]}{k+t} : |F \cap [k+t+s]| \ge k+t\right\} \subseteq \mathcal{F}$. 
For integers $m \ge r$, set families of the more general form 
$\left\{F \in \binom{[n]}{k} : |F \cap [m]| \ge r\right\}$ are widely considered to be extremal families in the literature of extremal set theory. 
Notable examples include stars, triangle families, and the Ahlswede--Khachatrian 
families (where one considers $\left\{F \in \binom{[n]}{k} : |F \cap [t+2i]| \ge t+i\right\}$). 
For more information on extremal constructions of this type, we refer the 
reader to a recent survey by Katona \cite{katonasurvey}. Furthermore, in 
\cite{li2013intersecting,wang2017intersecting,frankl2020erdHos}, 
several authors investigated the maximum size of an intersecting family 
containing a fixed family $\left\{F \in \binom{[n]}{k} : |F \cap [m]| \ge r\right\}$, 
a structural condition that Frankl \cite{frankl2020erdHos} referred to as 
a \emph{restricted universe}.

Building on these works, we develop a generalized framework that covers 
Conjecture \ref{frankl's-conj} as a special case (specifically, by setting 
the parameters to $k_1 = k_2+t-1$, $r = k_2+t-1$, and $m = k_2+t-1+s$ for $s\geq 0$). 
To state our main result, we need one more piece of notation. For any 
$A \subseteq [m]$ where $1\leq m<n$, we define
\[
A \times \binom{[m+1, n]}{x} = \left\{ A \cup B : B \subseteq [m+1, n],\ |B| = x \right\}.
\]
Our main results are as follows.

\begin{theorem}\label{mainthm}
    Let $n, k_1, k_2, m, r, t$ be integers satisfying $m \geq r \geq t$, 
    $k_1 > k_2 > k_1 - t$, $n \geq 2k_2 + t > k_1 + k_2 \geq r + k_2 \geq m + 1$, 
    and $2r \geq m + t$. 
    Let $\mathcal{F} \subseteq \binom{[n]}{k_1}$ and $\mathcal{G} \subseteq \binom{[n]}{k_2}$ 
    be cross-intersecting families. Suppose that 
    $\left\{F \in \binom{[n]}{k_1} : |F \cap [m]| \geq r\right\} \subseteq \mathcal{F}$ 
    and $\mathcal{F}$ is $t$-intersecting. Then
    \begin{equation*}
        |\mathcal{F}| + |\mathcal{G}| \leq 
        \left| \left\{ F \in \binom{[n]}{k_1} : |F \cap [m]| \geq r \right\} \right| 
        + \left| \left\{ G \in \binom{[n]}{k_2} : |G \cap [m]| \geq m + 1 - r \right\} \right|.
    \end{equation*}
    For $n \geq 2k_2 + t$, equality holds if and only if, up to permutations of the elements in $[n]$,    
    \begin{enumerate}[(a)]
        \item $\mathcal{F} = \left\{ F \in \binom{[n]}{k_1} : |F \cap [m]| \geq r \right\}$ 
        and $\mathcal{G} = \left\{ G \in \binom{[n]}{k_2} : |G \cap [m]| \geq m + 1 - r \right\}$.
        
        \item If $2(r-1) = t - 1 + m$ and $k_1 - k_2 = t - 1$, then besides the extremal 
        families in (a), we have $\mathcal{F} = \left\{ F \in \binom{[n]}{k_1} : |F \cap [m]| \geq r \right\} \cup \mathcal{H}$ 
        and $\mathcal{G} = \left\{ G \in \binom{[n]}{k_2} : |G \cap [m]| \geq m + 2 - r \right\} \cup \mathcal{K}$, 
        where $\mathcal{H}=\mathcal{P}\times\binom{[m+1,n]}{k_1-r+1}$ satisfies that $\mathcal{P}\subseteq\binom{[m]}{r-1}$ 
        is $t$-intersecting, and 
        $\mathcal{K} = \left( \binom{[m]}{m+1-r} \setminus \left\{ [m] \setminus P : P \in \mathcal{P} \right\} \right) \times \binom{[m+1, n]}{k_2 - (m + 1 - r)}$. 
    \end{enumerate}
\end{theorem}
We remark that in the Erd\H{o}s--Ko--Rado theorem, there exist extremal 
families other than stars when $n=2k$, a phenomenon that frequently occurs 
in extremal set theory. Note that in case (b) of Theorem \ref{mainthm}, 
by setting $k_2=s+2$, we obtain families $\mathcal{F}_0'$ and $\mathcal{G}_0'$ 
as follows: 
\begin{equation*}
    \mathcal{F}_0' = \left\{ F\in\binom{[n]}{s+t+1} : |F \cap [2s+t+1]|= s+t+1 \right\} \cup \mathcal{H},     
\end{equation*}
where $\mathcal{H}$ is such that $\{H \cap [2s+t+1]: |H \cap [2s+t+1]| = s+t, H\in \mathcal{H}\}$ 
is $t$-intersecting; and 
\begin{equation*}
    \mathcal{G}_0' = \left\{ G\in\binom{[n]}{s+2} : |G \cap [2s+t+1]|= s+2 \right\} \cup \left( \mathcal{K} \times \binom{[2s+t+2,n]}{1} \right),   
\end{equation*}
where $\mathcal{K}=\binom{[2s+t+1]}{s+1} \setminus \left\{ [2s+t+1] \setminus (H\cap [2s+t+1]) : H \in \mathcal{H} \right\}$.

One can verify that $\mathcal{F}_0'$ and $\mathcal{G}_0'$ also form extremal 
families achieving the upper bound in Conjecture \ref{frankl's-conj}. 
This construction illustrates why a complete characterization of all extremal 
families is likely to be complicated. We remark that the special case of 
this construction when $k_2 = 2$ was also presented in \cite{wu2023refined}.

\noindent \textbf{Proof Overview and New Techniques.} Our proof method builds on a 
framework combining the classical shifting operation of Erd\H{o}s--Ko--Rado 
\cite{erdos1961intersection} and the Ahlswede--Khachatrian generating set method 
\cite{ahlswede1997complete}. Within this framework, we go a step further by 
leveraging the normal property introduced in \cite{li2025non} to obtain a refined 
structural characterization of generating sets. 

The key step is to define critical pairs in normal families and derive 
generating-set-like properties. In previous literature, the generating set 
method could primarily handle intersecting families with relatively simple 
structures. Through a detailed analysis of critical pairs, we can construct 
new families from putative counterexamples to derive a contradiction whenever 
the original families deviate from our extremal families. Furthermore, 
we introduce operations that transform complex-structured families into highly 
normal ones. We believe that these critical structures within normal families 
are of independent interest and may find applications in related problems.

\noindent \textbf{Paper Organization.} The paper is organized as follows. 
Section~\ref{sec:prelim} extends the properties of generating sets from the 
classical settings of $t$-intersecting and cross-intersecting families to the 
mixed case. We then introduce the concept of normality, defines projected 
families, and derives generating-set-like properties, which are subsequently 
used to construct new families. In Section~3, we prove two technical lemmas 
demonstrating that achieving the maximum sum of sizes forces the two families 
to possess normal properties. Finally, the proof of Theorem~\ref{mainthm} 
is given in Section~\ref{proofofthm}.

\section{Preliminary Results}\label{sec:prelim}

In this section, we will collect and prove some basic properties of generating families and projected families in our setting, which will help construct new families in the proof of our main theorem.

\subsection{Generating Families}
Let $\mathcal{F}$ be a family of $k$-subsets of $[n]$. That is, $\mathcal{F}\subseteq\binom{[n]}{k}$. For $i,j\in[n]$ and $A\in \mathcal{F}$, define
$$
s_{i,j}(A)=\begin{cases}\begin{array}{cc}
    (A\setminus\{j\})\cup\{i\},&\text{if } j\in A,i\notin A, (A\setminus\{j\})\cup\{i\}\notin \mathcal{F};\\
    A,&\text{otherwise}
    \end{array}\end{cases}
$$
and $s_{i,j}(\mathcal{F})=\{s_{i,j}(A):A\in\mathcal{F}\}$. It is well known and easy to prove that
$s_{i,j}$ preserves the size $|\mathcal{F}|$ and the intersecting property. We apply the shifting operations to $\mathcal{F}$ for all $1\leq i<j\leq n$ until we obtain a family $\mathcal{F}'$ such that $s_{i,j}(\mathcal{F}')=\mathcal{F}'$ for all $1\leq i<j\leq n$; such an $\mathcal{F}'$ is called \emph{left-compressed}. Also, it is clear that $s_{i,j}(\mathcal{F}_1)$ and $s_{i,j}(\mathcal{F}_2)$ are cross-intersecting if and only if $\mathcal{F}_1$ and $\mathcal{F}_2$ are cross-intersecting.

It is clear that $s_{i,j}(\mathcal{F})$ and $s_{i,j}(\mathcal{G})$ are cross-intersecting with $T=\{A\in\binom{[n]}{k}:|A\cap [m]|\geq r\}\subseteq s_{i,j}(\mathcal{F})$ if and only if $\mathcal{F}$ and $\mathcal{G}$ are cross-intersecting with $T=\{A\in\binom{[n]}{k}:|A\cap [m]|\geq r\}\subseteq \mathcal{F}$. Therefore, for cross-intersecting families $\mathcal{F},\mathcal{G}$ satisfying $T \subseteq \mathcal{F}$, we may assume without loss of generality that $\mathcal{F},\mathcal{G}$ are left-compressed with $T \subseteq \mathcal{F}$.

In their seminal work \cite{ahlswede1997complete}, Ahlswede and Khachatrian introduced the following framework, developed the generating set method for $t$-intersecting families, and proved the complete EKR theorem. 

Let $\mathcal{F} \subseteq \binom{[n]}{k}$ be a left-compressed family. A family $g(\mathcal{F}) \subseteq \bigcup_{i \leq k} \binom{[n]}{i}$ is called a \emph{generating family} of $\mathcal{F}$ if 
\[
\left\{ B \in \binom{[n]}{k} : A \subseteq B,\ A \in g(\mathcal{F}) \right\} = \mathcal{F}.
\]
Let $G(\mathcal{F})$ denote the collection of all generating families of $\mathcal{F}$. Note that $\mathcal{F}$ itself is a generating family of $\mathcal{F}$, hence $G(\mathcal{F})\not=\emptyset$. For any $g(\mathcal{F}) \in G(\mathcal{F})$, let $g_{*}(\mathcal{F})$ be the set of minimal elements of $g(\mathcal{F})$ with respect to set inclusion. One can verify that $g_{*}(\mathcal{F}) \in G(\mathcal{F})$. We may therefore restrict our attention to 
$G_{*}(\mathcal{F}) = \{ g(\mathcal{F}) \in G(\mathcal{F}) : g(\mathcal{F}) = g_{*}(\mathcal{F}) \}$. The \emph{extent} of $g(\mathcal{F}) \in G_{*}(\mathcal{F})$ is defined as the maximum element appearing in $\{ \max(A) : A \in g(\mathcal{F}) \}$, and $\mathscr{D}(E) =\{F \in \mathcal{F} : F\cap   [\max{(E)}] =E\}$.

\begin{remark}
    While previous literature often uses the term ``generating set'', we 
    adopt the term ``generating family'' instead. Since $g(\mathcal{F})$ 
    is fundamentally a collection of sets rather than a set of elements, 
    this terminology more accurately reflects its mathematical structure.
\end{remark}

\begin{example}
    The following example helps illustrate these concepts.
    Consider the left-compressed family $\mathcal{F}=\left\{F\in\binom{[n]}{k}:1\in F\right\}$ where $n\geq k\geq 2$. Let $g(\mathcal{F})=\{\{1,2\},\{1,3\},\ldots,\\\{1,n\},\{1\}\}$. Clearly, $g(\mathcal{F})\in G(\mathcal{F})$ and $g_{*}(\mathcal{F})=\{\{1\}\}\in G_{*}(\mathcal{F})$. 
    Indeed, both $\{\{1,2\},\{1,3\},\\\ldots,\{1,n\}\}$ and $\{\{1\}\}$ are contained in $G_{*}(\mathcal{F})$; 
    the extent of $\{\{1\}\}$ is $1$ while the extent of $\{\{1,2\},\{1,3\},\ldots,\{1,n\}\}$ is $n$.
\end{example}

Let $\mathcal{F} \subseteq \binom{[n]}{k}$ be a left-compressed family, $g(\mathcal{F}) \in G_{*}(\mathcal{F})$ be a generating family of $\mathcal{F}$, and $e$ be the extent of $g(\mathcal{F})$. We define: 
\begin{align*}
g^{*}(\mathcal{F}) &= \left\{ A \in g(\mathcal{F}) : e \in A \right\}, \\
g_{u}^{*}(\mathcal{F}) &= \left\{ A \in g(\mathcal{F}) : |A| = u,\ e \in A \right\}, \\
g_{u}^{*}(\mathcal{F})' &= \left\{ A'=A \setminus \{\max(A)\} : A \in g_{u}^{*}(\mathcal{F}) \right\}.
\end{align*}
It is easy to check for different $E$ and $F$ in $g_{u}^{*}(\mathcal{F})$, $\mathscr{D}(E)$ and $\mathscr{D}(F)$ are different, hence $$\mathscr{D}(g_{u}^{*}(\mathcal{F})) = \bigcup_{E \in g_{u}^{*}(\mathcal{F})} \mathscr{D}(E).$$
The following key properties were established in \cite{ahlswede1997complete}; the authors of \cite{zhang2025conjecture} proved similar properties for 
cross-$t$-intersecting families. For completeness, we provide a proof of 
the following lemma in our setting.

\begin{lemma}{\em \cite{ahlswede1997complete,zhang2025conjecture}}\label{t-intersecting-cross-int-lemma}
    Let $\mathcal{F} \subseteq \binom{[n]}{k_1}$ and $\mathcal{G} \subseteq \binom{[n]}{k_2}$ be two left-compressed families satisfying the conditions in Theorem \ref{mainthm} and $|\mathcal{F}|+|\mathcal{G}|$ is maximized, with $g(\mathcal{F}) \in G_{*}(\mathcal{F})$ and $g(\mathcal{G}) \in G_{*}(\mathcal{G})$. Let $e = \max\{e_1,e_2\}$, where $e_1, e_2$ are the extents of $g(\mathcal{F})$ and $g(\mathcal{G})$, respectively. 
    \begin{enumerate}[(i)]
        \item For every $t < u < e$, $|g_{u}^{*}(\mathcal{F})|=|g_{e+1-u}^{*}(\mathcal{G})|$. Furthermore, if none of $g_{u}^{*}(\mathcal{F})$ and $g_{e+1-u}^{*}(\mathcal{G})$ is empty, $E_1 \in g_{u}^{*}(\mathcal{F})$ if and only if  $F_1 \in g_{e+1-u}^{*}(\mathcal{G})$ such that 
              $E_1 \cap F_1 = \{e\}$ and $E_1 \cup F_1 = [e]$. 
        \item Suppose $g_{u}^{*}(\mathcal{F}) \neq \emptyset$ for some integer $t < u < e$, and take $E_1 \in g_{u}^{*}(\mathcal{F})$. If there exists $E_2 \in g^{*}(\mathcal{F})$ such that $|E_1 \cap E_2| = t$, then $E_1 \cup E_2 = [e]$, and so $E_2\in g_{e+t-u}^{*}(\mathcal{F})$. 

        \item If $e>m$, $\mathscr{D}(g_{u}^{*}(\mathcal{F}))$ and $\mathscr{D}(g_{e+t-u}^{*}(\mathcal{F}))$ contain no elements of $\left\{ F : |F \cap [m]| \geq r \right\}$. 
    \end{enumerate}
\end{lemma}
\begin{proof}
    (i) Note that $G_{*}(\mathcal{F})$ contains the minimal elements with respect to set inclusion. Then $E_1 \in g_{u}^{*}(\mathcal{F})$ means that no proper subset of $E_1$ is in $g(\mathcal{F})$. Furthermore, the cross-intersecting property of $\mathcal{F}$ and $\mathcal{G}$ implies that $|F\cap E_1|\geq 1$ for any $F\in g(\mathcal{G})$. We claim that there must exist some $F_1 =\big([e] \setminus E_1'\big) \in g_{e+1-u}^{*}(\mathcal{G})$. If $F_1\notin g_{e+1-u}^{*}(\mathcal{G})$, we deduce that either there exists $F_0\subsetneq F_1$ with $F_0\in g(\mathcal{G})$, or there exists $E_3\in g(\mathcal{F})$ such that $E_3\cap F_1=\emptyset$ by the maximality of $|\mathcal{F}| + |\mathcal{G}|$. 
        In the former case, we know that $F_0\in g^{*}(\mathcal{G})$ since $|F_0\cap E_1|\geq 1$. Choose an integer $b\in F_1\setminus F_0$. Then there exists $F_3\in g(\mathcal{G})$ such that $F_3\subseteq s_{b,e}(F_0)$, yielding a contradiction because $F_3\cap E_1=\emptyset$. In the latter case, we have $E_3 \subsetneq E_1'$, which also yields a contradiction. Therefore, the claim is proved.

    Assume that $F\in g_{e+1-u}^{*}{(\mathcal{G})}$. Clearly there exists $E\in g(\mathcal{F})$ such that $E\cap F'=\emptyset$, otherwise there exists some subset of $F'$ in $g{(\mathcal{G})}$ by the maximality of $|\mathcal{F}| + |\mathcal{G}|$. Further, we deduce that $E\cap F=\{e\}$. We now show that $E \cup F= [e]$. Suppose not; then there exists some $b \in [e-1]$ such that $b \notin E \cup F$. Since $\mathcal{F}$ is left-compressed, there exists $B \in g(\mathcal{F})$ with $B \subseteq s_{b,e}(E)$. But then $B \cap F =\emptyset$, a contradiction.

    (ii) We now show that $E_1 \cup E_2 = [e]$. Suppose not; then there exists some $b \in [e-1]$ such that $b \notin E_1 \cup E_2$. Since $\mathcal{F}$ is left-compressed, there exists $B \in g(\mathcal{F})$ with $B \subseteq s_{b,e}(E_1)$. But then $|B \cap E_2| = t-1$, a contradiction.

    (iii) To establish this, note that $ \left\{ F : |F \cap [m]| \geq r \right\} \subseteq \mathcal{F}$. The condition $g(\mathcal{F}) \in G_{*}(\mathcal{F})$ implies that for any $r$-subset $R$ in $[m]$, there exists a subset $\overline{T} \subseteq R$ in $g(\mathcal{F})$ and the maximum element of $\overline{T}$ is at most $m$. By the definition of the generating family, $\mathcal{F}$ is a disjoint union of $\mathscr{D}(E)$, where $E\in g(\mathcal{F})$. Given that $e > m$, it is clear that $\overline{T} \notin \mathscr{D}(g_{u}^{*}(\mathcal{F}))$ or $\mathscr{D}(g_{e+t-u}^{*}(\mathcal{F}))$. Consequently, $\mathscr{D}(g_{u}^{*}(\mathcal{F}))$ and $\mathscr{D}(g_{e+t-u}^{*}(\mathcal{F}))$ contain no elements of $T$.
\end{proof}

Since $g_{u}^{*}(\mathcal{F}) \neq \emptyset$ for some $t<u<e$, we distinguish two cases according to whether $g_{e+t-u}^{*}(\mathcal{F})$ is empty or not:
    \begin{itemize}
        \item[(a)] $g_{u}^{*}(\mathcal{F}) \not= \emptyset$, $g_{e+t-u}^{*}(\mathcal{F})\not= \emptyset$, $g_{e+1-u}^{*}(\mathcal{G})\not= \emptyset$, and $g_{u+1-t}^{*}(\mathcal{G})\not= \emptyset$;
        \item[(b)] only $g_{u}^{*}(\mathcal{F})\not= \emptyset$ and $g_{e+1-u}^{*}(\mathcal{G})\not= \emptyset$.
    \end{itemize}
We construct families $\mathcal{F}_i$ and $\mathcal{G}_i$ for $i\in [4]$ in Table~\ref{tab:constructions} such that they are cross-intersecting, $\mathcal{F}_i$ is $t$-intersecting and $\left\{ F : |F \cap [m]| \geq r \right\} \subseteq \mathcal{F}_i$ when $e > m$ by Lemma \ref{t-intersecting-cross-int-lemma}.

\begin{table}[htbp]
\centering
\caption{Constructions of new cross-intersecting families}\label{tab:constructions}
\begin{tabular}{|c|c|c|}
\hline
Case & New families from $\mathcal{F}$ & New families from $\mathcal{G}$ \\
\hline
\multirow{2}{*}{(a)} 
& $\mathcal{F}_1 = \mathcal{F} \cup \mathscr{D}(g_{e+t-u}^{*}(\mathcal{F})') \setminus \mathscr{D}(g_{u}^{*}(\mathcal{F}))$ 
& $\mathcal{G}_1 = \mathcal{G} \cup \mathscr{D}(g_{e+1-u}^{*}(\mathcal{G})') \setminus \mathscr{D}(g_{u+1-t}^{*}(\mathcal{G}))$ \\
\cline{2-3}
& $\mathcal{F}_2 = \mathcal{F} \cup \mathscr{D}(g_{u}^{*}(\mathcal{F})') \setminus \mathscr{D}(g_{e+t-u}^{*}(\mathcal{F}))$ 
& $\mathcal{G}_2 = \mathcal{G} \cup \mathscr{D}(g_{u+1-t}^{*}(\mathcal{G})') \setminus \mathscr{D}(g_{e+1-u}^{*}(\mathcal{G}))$\\
\hline
\multirow{2}{*}{(b)} 
& $\mathcal{F}_3 = \mathcal{F} \setminus \mathscr{D}(g_{u}^{*}(\mathcal{F}))$ 
& $\mathcal{G}_3 = \mathcal{G} \cup \mathscr{D}(g_{e+1-u}^{*}(\mathcal{G})')$ \\
\cline{2-3}
& $\mathcal{F}_4 = \mathcal{F} \cup \mathscr{D}(g_{u}^{*}(\mathcal{F})')$ 
& $\mathcal{G}_4 = \mathcal{G} \setminus \mathscr{D}(g_{e+1-u}^{*}(\mathcal{G}))$ \\
\hline
\end{tabular}
\end{table}

\subsection{Normality and Projected Families}

This section is devoted to establishing normality results for cross-intersecting families. After introducing the definition of normal families, we show some basic properties, and prove that $\mathcal{F}$ and $\mathcal{G}$ are $[m+1,n]$-normal (Lemma~\ref{[m+1-n]-unifrom-lemma}). Building on this, we assume $[\ell+1,m]$-normality to establish generating-set-like properties (Lemma~\ref{key-claim-frankl}). These properties are then used to construct  new families and used to prove that $[1,m]$-normal families $\mathcal{F}$ and $\mathcal{G}$ are one of the largest (Lemma~\ref{[1,m]-normal}). Lemmas \ref{[m+1-n]-unifrom-lemma} and \ref{[1,m]-normal} will be proved in Section 3.

For $S\subseteq[n]$ and $A,B\in \binom{[n]}{k}$, we say that $A$ and $B$ are $S$-conjugate if $A\setminus S=B\setminus S$. This is an equivalence relation on $\binom{[n]}{k}$. For left-compressed family $\mathcal{F}\subseteq\binom{[n]}{k}$, we say $\mathcal{F}$ is \emph{S-normal} if $\mathcal{F}$ is a union of some of $S$-conjugate classes.  
More precisely, let $u(\mathcal{F})_{S}$ collect all $F\cap S$, where $F\in \mathcal{F}$. If $\mathcal{F}$ is \emph{S-normal}, then
\[
\mathcal{F} = \bigcup_{ R\in u(\mathcal{F})_{[n]\setminus S}}R\times\binom{S}{k - |R|}.
\]

\begin{lemma}\label{[m+1-n]-unifrom-lemma}
 Let $\mathcal{F} \subseteq \binom{[n]}{k_1}$ and $\mathcal{G} \subseteq \binom{[n]}{k_2}$ be two left-compressed families satisfying the conditions in Theorem \ref{mainthm} and $|\mathcal{F}|+|\mathcal{G}|$ is maximized. Then $\mathcal{F}$ and $\mathcal{G}$ are $[m+1,n]$-normal.
\end{lemma}

\begin{lemma}\label{[1,m]-normal}
 There exist two left-compressed $[1,m]$-normal families $\mathcal{F} \subseteq \binom{[n]}{k_1}$ and $\mathcal{G} \subseteq \binom{[n]}{k_2}$ satisfying the conditions in Theorem \ref{mainthm} and $|\mathcal{F}|+|\mathcal{G}|$ is maximized.  
\end{lemma}

Recall that $e = \max\{e_1,e_2\}$, where $e_1, e_2$ are the extents of $g(\mathcal{F})$ and $g(\mathcal{G})$, respectively. Then by the decomposing properties of generating sets in \cite{ahlswede1997complete,zhang2025conjecture} and Lemma \ref{t-intersecting-cross-int-lemma},  $\mathcal{F}$ and $\mathcal{G}$ are $[e+1,n]$-normal. Lemma~\ref{[m+1-n]-unifrom-lemma} further shows that $e \leq m$.

Among all projected families $u(\mathcal{F})_{S}$, we will frequently use $u(\mathcal{F})_{[m]}$ in the rest of the paper, so we omit $[m]$ and denote
\[
u(\mathcal{F}) = \{ F \cap [m] : F \in \mathcal{F} \}
\quad \text{and} \quad
u(\mathcal{G}) = \{ G \cap [m] : G \in \mathcal{G} \}.
\]
Then by Lemma~\ref{[m+1-n]-unifrom-lemma}, the families $\mathcal{F}$ and $\mathcal{G}$ are $[m+1,n]$-normal, so they admit the decompositions
\[
\mathcal{F} = \bigcup_{A \in u(\mathcal{F})} A \times \binom{[m+1, n]}{k_1 - |A|}
\quad \text{and} \quad
\mathcal{G} = \bigcup_{B \in u(\mathcal{G})} B \times \binom{[m+1, n]}{k_2 - |B|}.
\]

For convenience, we write $A \times \binom{[\ell+1, m]}{x} \subseteq u(\mathcal{F})$ if every set in this collection belongs to $u(\mathcal{F})$. Conversely, $A \times \binom{[\ell+1, m]}{x}\cap u(\mathcal{F})=\emptyset$ means that every set in this collection is not
contained in $u(\mathcal{F})$. We now show that the above notation is crucial for understanding the normality inside $[1,m]$. Suppose that $\mathcal{F}$ is $[\ell+1, m]$-normal but not $[\ell, m]$-normal for some positive integer $\ell< m$. Then there exist a set $E_1 \subseteq [\ell]$ with $\ell \in E_1$, and an integer $0 \le x \le m-\ell-1$:
\[
E_1 \times \binom{[\ell+1, m]}{x} \subseteq u(\mathcal{F})
\text{ and } 
E_1' \times \binom{[\ell+1, m]}{x+1} \cap u(\mathcal{F}) = \emptyset, \text{ where } E_1'=E_1\setminus\{\ell\}.
\]
The existence of such a $E_1$ is the key to constructing maximized families. For a fixed $\ell$, we call a pair $(E_1,x)$ critical in $\mathcal{F}$ if they satisfy the two conditions above.
Let $u_{a}^{*}(\mathcal{F}) = \{E_1\mid |E_1|=a, \text{$(E_1,x)$  is critical in $\mathcal{F}$ for some nonnegative integer $x$}\}$. For any $E_1\in u_{a}^{*}(\mathcal{F})$, let $I(E_1)$ collect all $x$ such that $(E_1,x)$ is critical in $\mathcal{F}$. One can similarly define critical pairs $(F_1,y)$, $\I(F_1)$ and $u_{a}^{*}(\mathcal{G})$ in $\mathcal{G}$ for the same $\ell$.

We now use the normality property to refine the structure of $u(\mathcal{F})$ and establish the following key lemma.

\begin{lemma}\label{key-claim-frankl}
    Let $\mathcal{F} \subseteq \binom{[n]}{k_1}$ and $\mathcal{G} \subseteq \binom{[n]}{k_2}$ be two left-compressed families satisfying the conditions in Theorem \ref{mainthm} and $|\mathcal{F}|+|\mathcal{G}|$ is maximized. 
    Suppose further that $\mathcal{F}$ is $[\ell+1, m]$-normal but not $[\ell, m]$-normal for some positive integer $\ell< m$. Then the following hold:
    
    \begin{enumerate}[(i)]

        \item $u(\mathcal{F})$ and $u(\mathcal{G})$ are cross intersecting, and $u(\mathcal{F})$ is $t$-intersecting.
        
        \item For each $E_1\in u_{a}^{*}(\mathcal{F})$, $I(E_1)$ is an interval. 

\item  $\mathcal{G}$ is $[\ell+1, m]$-normal. 

     \item  $E_1\in u_{a}^{*}(\mathcal{F})$ if and only if $F_1\in u_{\ell-a+1}^{*}(\mathcal{G})$. 
Furthermore,  $I(F_1)$ is the interval consisting exactly of all $y=m-\ell-1-x$ with $x \in I(E_1)$. 

    \end{enumerate}
\end{lemma}
\begin{proof}[\bf Proof of Lemma \ref{key-claim-frankl}.]
    We first prove part (i).  It follows from $e\le m$ by Lemma~\ref{[m+1-n]-unifrom-lemma} and the fact that the generating families keep the intersecting properties.

    For part (ii), it is easy to see that  $E_1 \times \binom{[\ell+1, m]}{x} \subseteq u(\mathcal{F})$ implies that $E_1 \times \binom{[\ell+1, m]}{x+1} \subseteq u(\mathcal{F})$ 
    by the maximality of $|\mathcal{F}|+|\mathcal{G}|$. 
    By the left-compressed property, $E_1' \times \binom{[\ell+1, m]}{x+1} \subseteq u(\mathcal{F})$ implies 
    $E_1\times \binom{[\ell+1, m]}{x} \subseteq u(\mathcal{F})$. 
    Combining these two facts, we deduce the existence of $I(E_1)$ and it is an interval.
    
    For part (iii), let $B=B_1\cup B_2\in u(\mathcal{G})$ where $B_1\subseteq [\ell]$ and $B_2\subseteq [\ell+1,m]$. Since $\mathcal{F}$ is $[\ell+1,m]$-normal, we have that
$u(\mathcal{F})$ consists of $A_1 \times \binom{[\ell+1, m]}{x}$ for some $A_1$ and $x$. By the maximality of $|\mathcal{F}|+|\mathcal{G}|$, if $A_1 \cap B_1 \not= \emptyset$, then $B_1\cup B_3\in u(\mathcal{G})$ for any $|B_2|=|B_3|$; otherwise $|B_2|+x>m-\ell$, we again obtain $B_1\cup B_3\in u(\mathcal{G})$ for any $|B_2|=|B_3|$.

    For part (iv), if $(E_1,x)$ is critical in $\mathcal{F}$, then $F_1' \times \binom{[\ell+1, m]}{m - \ell - x} \cap u(\mathcal{G}) = \emptyset$. Assume that $F_1 \times \binom{[\ell+1, m]}{m - \ell - x - 1} \cap u(\mathcal{G}) = \emptyset$. By the maximality of $|\mathcal{F}|+|\mathcal{G}|$, we deduce that there exists $A\in u(\mathcal{F})$ such that $A\cap\left( F_1 \times \binom{[\ell+1, m]}{m - \ell - x - 1}\right)=\emptyset$. Let $A=A_1\cup A_2$ where $A_1\subseteq [\ell]$ and $A_2\subseteq [\ell+1,m]$. Then $A_1\subseteq E_1'$ and $|A_2|\leq x+1$, which contradicts $E_1' \times \binom{[\ell+1, m]}{x+1} \cap u(\mathcal{F}) = \emptyset$. This contradiction shows $(F_1,m-\ell-x-1)$ is critical in $\mathcal{G}$.

    On the other direction, $F_1' \times \binom{[\ell+1, m]}{m-\ell-x} \cap u(\mathcal{G}) = \emptyset$ means that there exists $A\in u(\mathcal{F})$ such that $A\cap\left( F_1' \times \binom{[\ell+1, m]}{m-\ell-x}\right)=\emptyset$. Let $A=A_1\cup A_2$ where $A_1\subseteq [\ell]$ and $A_2\subseteq [\ell+1,m]$. Then $A_1\subseteq E_1$ and $|A_2|\leq x$. The maximality of $|\mathcal{F}|+|\mathcal{G}|$ implies that $E_1\times \binom{[\ell+1,m]}{x}\subseteq u(\mathcal{F})$. The cross-intersecting property of $\mathcal{F}$ and $\mathcal{G}$ means $E_1' \times \binom{[\ell+1, m]}{x +1} \cap u(\mathcal{F}) = \emptyset$. Furthermore, from the above proof, we have $|u_{a}^{*}(\mathcal{F})| = |u_{\ell-a+1}^{*}(\mathcal{G})|$ and $|I(E_1)| = |I(F_1)|$.
\end{proof}

As a consequence, we have the following lemma.

\begin{lemma}
    Let $\mathcal{F} \subseteq \binom{[n]}{k_1}$ and $\mathcal{G} \subseteq \binom{[n]}{k_2}$ be two left-compressed families satisfying the conditions in Theorem \ref{mainthm} and $|\mathcal{F}|+|\mathcal{G}|$ is maximized. 
    Then $\mathcal{F}$ is $[\ell+1, m]$-normal but not $[\ell, m]$-normal if and only if $\mathcal{G}$ is $[\ell+1, m]$-normal but not $[\ell, m]$-normal.

\end{lemma}

Now, we may assume that $\mathcal{F}$ and $\mathcal{G}$ are $[\ell+1, m]$-normal but not $[\ell, m]$-normal for some positive integer $1\le \ell\le m$. Let $E_1\in u_{a}^{*}(\mathcal{F})$ for some positive integer $1\leq a\leq \ell$. If there are $E_2$'s in $u_{\ell-a+t}^{*}(\mathcal{F})$ satisfying $|E_1 \cap E_2| = t$ and $E_1 \cup E_2 = [\ell]$, we call $(E_1,E_2)$ is a good pair. 

\begin{lemma}\label{key-claim-E2}
There exist two left-compressed families $\mathcal{F} \subseteq \binom{[n]}{k_1}$ and $\mathcal{G} \subseteq \binom{[n]}{k_2}$ and an integer $1\le \ell\le m$ satisfying the conditions in Theorem \ref{mainthm} such that
 $|\mathcal{F}|+|\mathcal{G}|$ is maximized and for any $E_1 \in u_{a}^{*}(\mathcal{F})$ with an interval $I(E_1)$, if there exists a good pair $(E_1,E_2)$, then $I(E_2)$ is an interval that collects all $m-1-\ell-x$, where $x \in I(E_1)$. 
\end{lemma}

Before we prove Lemma \ref{key-claim-E2}, we first show the following two claims, which will help us to construct new cross-intersecting families.

\begin{claim} Suppose there exist a set $A$ and an integer $a$ such that  
    $A \times \binom{[\ell+1,m]}{a} \subseteq u(\mathcal{G})$, 
and such that the two families  
$A \times \binom{[\ell+1,m]}{a}, E_1 \times \binom{[\ell+1,m]}{x}$
are intersecting, while  
$A \times \binom{[\ell+1,m]}{a}, E_1' \times \binom{[\ell+1,m]}{x+1}$
are not intersecting for some $x\in I(E_1)$.  
Then $A=[\ell]\setminus E_1'$ and $a \le m - \ell - (x+1)$.
\end{claim}

\noindent{\bf Proof of Claim 1.} We deduce that $\ell \in A$, $|A \cap E_1| = 1$, and $a \leq m-\ell-x-1$. Now assume for a contradiction that $A \cup E_1 \subsetneq [\ell]$. Then there exists an integer $c \in [\ell] \setminus \{A \cup E_1\}$. By the left-compressed property of $\mathcal{F}$, we have
    $s_{c,\ell}\Bigl(E_1 \times \binom{[\ell+1,m]}{x}\Bigr) \subseteq u(\mathcal{F})$.
    But this implies that the families
    $A \times \binom{[\ell+1,m]}{a}, s_{c,\ell}\Bigl(E_1 \times \binom{[\ell+1,m]}{x}\Bigr)$ are not intersecting, which contradicts our earlier assumption. $\hfill\blacksquare$
    
\begin{claim}
    
Suppose there exist a set $B$ and an integer $b$ such that  
$B \times \binom{[\ell+1,m]}{b} \subseteq u(\mathcal{F})$, 
and such that the two families  
$B \times \binom{[\ell+1,m]}{b}, E_1 \times \binom{[\ell+1,m]}{x}$
are $t$-intersecting, while  
$B \times \binom{[\ell+1,m]}{b}, E_1' \times \binom{[\ell+1,m]}{x+1}$
are not $t$-intersecting for some $x\in I(E_1)$.     
Then $B$ must satisfy  
$\ell \in B,|B \cap E_1| = t, B \cup E_1 = [\ell],$
and furthermore  
$b \le m - \ell - (x+1)$.
\end{claim} 

\noindent{\bf Proof of Claim 2.} We conclude that $\ell \in B$, $|B \cap E_1| = t$, and $b \leq m-\ell-x-1$. Now, suppose for contradiction that $B \cup E_1 \subsetneq [\ell]$. Then there exists an integer $c \in [\ell] \setminus \{B \cup E_1\}$. By the left-compressed property of $\mathcal{F}$, we have
    $s_{c,\ell}\Bigl(E_1 \times \binom{[\ell+1,m]}{x}\Bigr) \subseteq u(\mathcal{F})$.
    But this implies that the families
    $B \times \binom{[\ell+1,m]}{b}, s_{c,\ell}\Bigl(E_1 \times \binom{[\ell+1,m]}{x}\Bigr)$ are not $t$-intersecting, which contradicts our earlier assumption.  $\hfill\blacksquare$

Now we are ready to prove Lemma \ref{key-claim-E2}. 

\begin{proof}[Proof of Lemma \ref{key-claim-E2}]

We distinguish two cases according to whether or not there exists a set $W$ and an integer $w$ such that  
$W \times \binom{[\ell+1,m]}{w} \subseteq u(\mathcal{F})$, 
and such that the two families  
$W \times \binom{[\ell+1,m]}{w}, E_1\times \binom{[\ell+1,m]}{x_2}$
are $t$-intersecting, while  
$W \times \binom{[\ell+1,m]}{w}, E_1' \times \binom{[\ell+1,m]}{x_2+1}$
are not $t$-intersecting. Such condition is denoted as critical condition holds for $E_1$.

\begin{case}
Critical condition does not hold for $E_1$. 
\end{case} 

In this case, for any $W\times\binom{[\ell+1,m]}{w}\subseteq u(\mathcal{F})$, we have $W\times\binom{[\ell+1,m]}{w},E_1\times \binom{[\ell+1,m]}{x_2}$ are $t$-intersecting and $W\times\binom{[\ell+1,m]}{w},E_1'\times \binom{[\ell+1,m]}{x_2+1}$ are also $t$-intersecting. Assume that $\bar{x}$ is the largest integer such that for any $W\times\binom{[\ell+1,m]}{w}\subseteq u(\mathcal{F})$, we have $W\times\binom{[\ell+1,m]}{w},E_1\times \binom{[\ell+1,m]}{\bar{x}}$ are $t$-intersecting and $W\times\binom{[\ell+1,m]}{w},E_1'\times \binom{[\ell+1,m]}{\bar{x}+1}$ are not $t$-intersecting. By Claim 2, there exists $E_2\times\binom{[\ell+1,m]}{m-\ell-\bar{x}-1}\subseteq u(\mathcal{F})$ and $m-\ell-\bar{x}-1$ is the smallest non-negative integer. 

\begin{subcase}
    $x_2>\bar{x}\geq x_1$.
\end{subcase}

In this subcase, we have the following structure:
\begin{align*}
&E_1\times \binom{[\ell+1,m]}{x_1}\subseteq u(\mathcal{F}),\ldots, E_1\times \binom{[\ell+1,m]}{\bar{x}}\subseteq u(\mathcal{F})\\
&E_2\times \binom{[\ell+1,m]}{m-\ell-x_1-1}\subseteq u(\mathcal{F}),\ldots, E_2\times \binom{[\ell+1,m]}{m-\ell-\bar{x}-1}\subseteq u(\mathcal{F}).
\end{align*}
Now we can consider the following two constructions. Construction 1: $\mathcal{F}_1=\mathcal{F}\setminus \cup_{x\in [x_1,\bar{x}]}E_1\times\binom{[\ell+1,m]}{x}\times\binom{[m+1,n]}{k_1-a-x}$ and $\mathcal{G}_1=\mathcal{G}\cup_{x\in [x_1,\bar{x}]}F_1'\times\binom{[\ell+1,m]}{m-\ell-x}\times\binom{[m+1,n]}{k_2+a+x-m}$. These families are cross-intersecting by Claim 1 and the maximum of $|\mathcal{F}|+|\mathcal{G}|$ implies  
\begin{equation}\label{lem-ineq-3}
    \sum_{x\in [x_1,\bar{x}]}\binom{m-\ell}{x}\binom{n-m}{k_1-a-x}\geq \sum_{x\in [x_1,\bar{x}]}\binom{m-\ell}{x}\binom{n-m}{k_2+a+x-m}.
\end{equation}
Construction 2: $\mathcal{F}_2=\mathcal{F}\setminus \cup_{x\in [x_1,\bar{x}]}E_2\times\binom{[\ell+1,m]}{m-\ell-x-1}\times\binom{[m+1,n]}{k_1+a+x-m-(t-1)}$ and $\mathcal{G}_2=\mathcal{G}\cup\cup_{x\in [x_1,\bar{x}]}F_2'\times\binom{[\ell+1,m]}{x+1}\times\binom{[m+1,n]}{k_2-(a+x)+t-1}$. These families are cross-intersecting by Claim 1 and Claim 2, and the maximum of $|\mathcal{F}|+|\mathcal{G}|$ implies 
\begin{equation}\label{lem-ineq-4}
    \sum_{x\in [x_1,\bar{x}]}\binom{m-\ell}{x+1}\binom{n-m}{k_1+a+x-m-(t-1)}\geq \sum_{x\in [x_1,\bar{x}]}\binom{m-\ell}{x+1}\binom{n-m}{k_2-(a+x)+t-1}.
\end{equation}

Combining inequalities (\ref{lem-ineq-3}) and (\ref{lem-ineq-4}), we derive a contradiction if $|\bar{x}-x_1|\geq 1$. And if $\bar{x}=x_1$, the two inequalities hold if and only if $2(a+\bar{x})=k_1-k_2+m$. The detailed verification of the above facts, involving explicit estimates, is provided in Appendix (Lemma~\ref{fact-4}). In this equality case, we just delete  $E_2\times\binom{[\ell+1,m]}{m-\ell-1-x_1}$. If $E_2$ is not critical for any $x$, then $E_2$ is not critical for $[\ell+1,m]$-normality of $\mathcal{F}$. Otherwise, we have $(E_2, m-\ell-1-x_1)$ is a critical pair, that is the following subcase. 

\begin{subcase}
    $\bar{x}<x_1$. 
\end{subcase}

In this subcase, for any $W\times\binom{[\ell+1,m]}{w}\subseteq u(\mathcal{F})$, $W\times\binom{[\ell+1,m]}{w},E_1\times \binom{[\ell+1,m]}{x_1}$ are $t$-intersecting, and $W\times\binom{[\ell+1,m]}{w},E_1'\times \binom{[\ell+1,m]}{x_1+1}$ are also $t$-intersecting. Now we can consider the following two constructions. Construction 1: $\mathcal{F}_1=\mathcal{F}\setminus \cup_{x\in I(E_1)}E_1\times\binom{[\ell+1,m]}{x}\times\binom{[m+1,n]}{k_1-a-x}$ and $\mathcal{G}_1=\mathcal{G}\cup_{x\in I(E_1)}F_1'\times\binom{[\ell+1,m]}{m-\ell-x}\times\binom{[m+1,n]}{k_2+a+x-m}$. Those are cross-intersecting by Claim 1 and the maximum of $|\mathcal{F}|+|\mathcal{G}|$ implies  
\begin{equation}\label{lem-ineq-5}
    \sum_{x\in I(E_1)}\binom{m-\ell}{x}\binom{n-m}{k_1-a-x}\geq \sum_{x\in I(E_1)}\binom{m-\ell}{x}\binom{n-m}{k_2+a+x-m}.
\end{equation}
Construction 2: $\mathcal{F}_1=\mathcal{F} \cup_{x\in I(E_1)}E_1'\times\binom{[\ell+1,m]}{x+1}\times\binom{[m+1,n]}{k_1-a-x}$ and $\mathcal{G}_1=\mathcal{G}\setminus\cup_{x\in I(E_1)}F_1\times\binom{[\ell+1,m]}{m-\ell-x-1}\times\binom{[m+1,n]}{k_2+a+x-m}$. Those are cross-intersecting by Claim 1 and the maximum of $|\mathcal{F}|+|\mathcal{G}|$ implies 
\begin{equation}\label{lem-ineq-6}
    \sum_{x\in I(E_1)}\binom{m-\ell}{x+1}\binom{n-m}{k_2+a+x-m}\geq \sum_{x\in I(E_1)}\binom{m-\ell}{x+1}\binom{n-m}{k_1-a-x}.
\end{equation}
Combining inequalities (\ref{lem-ineq-5}) and (\ref{lem-ineq-6}), we derive a contradiction if $|I(E_1)|\geq 2$ by Lemma~\ref{key-binomial-ineq} and Lemma~\ref{fact-2}. And if $|I(E_1)|=1$, the equality holds if and only if $2(a+x)=k_1-k_2+m$. We just delete $E_1\times\binom{[\ell+1,m]}{x}$ such that $E_1$ is not critical for $[\ell+1,m]$-normality of $\mathcal{F}$.

\begin{case}
Critical condition holds for $E_1$. 
\end{case}

By Claim 2, we derive that $E_2\times \binom{[\ell+1,m]}{m-\ell-x_2-1}\subseteq u(\mathcal{F})$ for such $E_2$ with $\ell \in E_2,|E_2 \cap E_1| = t$ and $E_2 \cup E_1 = [\ell]$. Then the maximum of $|\mathcal{F}|+|\mathcal{G}|$ and the $t$-intersecting property of $\mathcal{F}$ imply
\begin{align*}
&E_2\times \binom{[\ell+1,m]}{m-\ell-x_1-1}\subseteq u(\mathcal{F}),\ldots, E_2\times \binom{[\ell+1,m]}{m-\ell-x_2-1}\subseteq u(\mathcal{F})  \\
&E_2'\times \binom{[\ell+1,m]}{m-\ell-x_1}\cap u(\mathcal{F}) = \emptyset,\ldots,E_2'\times \binom{[\ell+1,m]}{m-\ell-x_2}\cap u(\mathcal{F}) = \emptyset.
\end{align*}
Clearly, $|I(E_2)| \geq |I(E_1)|$. Now, we start from $E_2$. If $|I(E_1)| = |I(E_2)|$, then the desired conclusion follows. Otherwise, we can apply an argument similar to that in Case~1 to deduce that the lengths of $I(E_1)$ and $I(E_2)$ are equal. Since there are finitely many critical sets and the chosen complexity decreases at each step, the process terminates. Hence for any $E_1 \in u_{a}^{*}(\mathcal{F})$ associated with an interval $I(E_1)$, if there exists a good pair $(E_1, E_2)$, then $I(E_2)$ is an interval consisting of all elements of the form $m - 1 - \ell - x$, where $x \in I(E_1)$. This completes the proof of Lemma~\ref{key-claim-E2}.
\end{proof}

For every $E_1$ of size $a$, depending on whether there exists an $E_2$ such that $(E_1,E_2)$ is a good pair, we distinguish two cases:
    \begin{itemize}
        \item[(a)] There are no good pairs for every $E_1$;
        \item[(b)] There exists a good pair.
    \end{itemize}
For $E_1 \in u_a^*(\mathcal{F})$, we can find all $E_2$ such that $(E_1,E_2)$ is a good pair. Then for every $E_2$, we can find all $E_3$ such that $(E_3,E_2)$ is a good pair. We continue this procedure, and let $G(E_1)$ and $G(E_2)$ be all critical sets with the same size of $E_1$ and $E_2$, respectively. Now, we introduce the following notation for the two cases:
\[
\begin{cases}
\mathscr{D}(I(E_1)) = \left\{\bigcup_{x \in I(E_1)} E_1\times\binom{[\ell+1,m]}{x}\times\binom{[m+1,n]}{k_1-a-x}\right\},\\
\mathscr{D}(I(E_1)') = \left\{\bigcup_{x \in I(E_1)} E_1'\times\binom{[\ell+1,m]}{x+1}\times\binom{[m+1,n]}{k_1-a-x}\right\},\\
\mathscr{G}(I(E_1)) =  \left\{\bigcup_{E_1\in G(E_1)}\bigcup_{x \in I(E_1)} E_1 \times \binom{[\ell+1,m]}{x} \times \binom{[m+1,n]}{k_1-a-x} \right\},\\ \mathscr{G}(I(E_1)') =  \left\{\bigcup_{E_1\in G(E_1)}\bigcup_{x \in I(E_1)} E'_1 \times \binom{[\ell+1,m]}{x+1} \times \binom{[m+1,n]}{k_1-a-x} \right\}.
\end{cases}
\]
We use the analogous notation for $\mathcal{G}$.

We construct families $\widetilde{\mathcal{F}}_i$ and $\widetilde{\mathcal{G}}_i$ for $i\in [5]$ in Table~\ref{tab:constructions2} such that they are cross-intersecting, $\widetilde{\mathcal{F}}_i$ is $t$-intersecting and $\left\{ F : |F \cap [m]| \geq r \right\} \subseteq \widetilde{\mathcal{F}}_i$  by Lemma \ref{key-claim-frankl} and the above two claims.

\begin{table}[htbp]
\centering
\renewcommand{\arraystretch}{1.2} 
\caption{Constructions of new cross-intersecting families}\label{tab:constructions2}
\begin{tabular}{|c|c|c|}
\hline
Case & New families from $\mathcal{F}$ & New families from $\mathcal{G}$ \\
\hline
\multirow{2}{*}{(a)} 
& $\widetilde{\mathcal{F}}_1 = \mathcal{F} \setminus \mathscr{D}(I(E_1))$ 
& $\widetilde{\mathcal{G}}_1 = \mathcal{G} \cup \mathscr{D}(I(F_1)')$  \\
\cline{2-3}
& $\widetilde{\mathcal{F}}_2 = \mathcal{F} \cup \mathscr{D}(I(E_1)')$ 
& $\widetilde{\mathcal{G}}_2 = \mathcal{G} \setminus \mathscr{D}(I(F_1))$\\
\hline
\multirow{3}{*}{(b)} 
& $\widetilde{\mathcal{F}}_3 = \mathcal{F} \setminus \mathscr{G}(I(E_1)) \cup \mathscr{G}(I(E_2)')$ 
& $\widetilde{\mathcal{G}}_3 = \mathcal{G} \setminus \mathscr{G}(I(F_2)) \cup \mathscr{G}(I(F_1)')$ \\
\cline{2-3}
& $\widetilde{\mathcal{F}}_4 = \mathcal{F} \setminus \mathscr{G}(I(E_2)) \cup \mathscr{G}(I(E_1)')$ 
& $\widetilde{\mathcal{G}}_4 = \mathcal{G} \setminus \mathscr{G}(I(F_1)) \cup \mathscr{G}(I(F_2)')$ \\
\cline{2-3}
& $\widetilde{\mathcal{F}}_5 = \mathcal{F} \setminus \mathscr{G}(I(E_1))$ 
&  $\widetilde{\mathcal{G}}_5 = \mathcal{G} \cup \mathscr{G}(I(E_2)')$ \\
\hline
\end{tabular}
\end{table}

\subsection{A Combinatorial Inequality}

In this subsection, we prove one more technical lemma. Let $I = \{a, a+1, \ldots, b\}$ 
be an interval of integers contained in $[n]$ with $|I| \geq 2$. Let $(s_x)_{x \in I}$ be a sequence of 
real numbers such that $s_a > 0$ and the sign changes at most once; that is, there 
exists an integer $k \in I$ satisfying $s_x \geq 0$ for all $x \leq k$ and 
$s_x \leq 0$ for all $x > k$. In our proof, we make use of the following inequalities 
involving binomial coefficients.

\begin{lemma}\label{key-binomial-ineq}
    If $\sum_{x\in I}\binom{n}{x}s_x \geq 0$, then $\sum_{x\in I}\binom{n}{x+1}s_x > 0$. 
\end{lemma}

\begin{proof}
    Let $I_1 = \{x \in I : s_x \geq 0\}$ and $I_2 = \{y \in I : s_y < 0\}$. 
    If $I_2 = \emptyset$, then $s_x \geq 0$ for all $x \in I$. Since $s_a > 0$, 
    it trivially follows that $\sum_{x \in I} \binom{n}{x+1} s_x > 0$, and we are done. 
    
    Thus, we may assume $I_2 \neq \emptyset$. We prove the statement by contradiction. 
    Suppose that 
    \begin{equation*}
        \sum_{x\in I_1}\binom{n}{x}|s_x| \geq \sum_{y\in I_2}\binom{n}{y}|s_y|
    \end{equation*}
    but
    \begin{equation*}
        \sum_{y\in I_2}\binom{n}{y+1}|s_y| \geq \sum_{x\in I_1}\binom{n}{x+1}|s_x|. 
    \end{equation*}
    Multiplying the above two inequalities, we obtain
    \begin{equation*}
        \sum_{x\in I_1, y\in I_2}\binom{n}{x}\binom{n}{y+1}|s_x||s_y| \geq \sum_{x\in I_1, y\in I_2}\binom{n}{y}\binom{n}{x+1}|s_x||s_y|.
    \end{equation*}
   
    Notice that the function $f(z) = \frac{\binom{n}{z}}{\binom{n}{z+1}} = \frac{z+1}{n-z}$ 
    is strictly increasing with respect to $z$. The condition that the sign of $s_x$ 
    changes at most once implies that $\max I_1 < \min I_2$. Therefore, we have 
    $x < y$ for all $x \in I_1$ and $y \in I_2$, yielding $f(x) < f(y)$, or equivalently, 
    $\binom{n}{x}\binom{n}{y+1} < \binom{n}{y}\binom{n}{x+1}$.
    
    Since $s_a > 0$ and $I_2 \neq \emptyset$, there exists at least one pair $(x,y)$ 
    such that $|s_x||s_y| > 0$. This forces the strict inequality to hold for the 
    entire sum, which contradicts our earlier derived inequality.
\end{proof}

\section{Proofs  of Lemmas on Normality}\label{sec:normal}

In this section, we prove the two technical lemmas. 

\subsection{Families $\mathcal{F}$ and $\mathcal{G}$ are $[m+1,n]$-normal}

\begin{proof}[Proof of Lemma \ref{[m+1-n]-unifrom-lemma}]
Recall that $T=\left\{F\in\binom{n}{k_1}:| F\cap [m]|\geq r\right\}\subseteq \mathcal{F}$, with $g(\mathcal{F})\in G_{*}(\mathcal{F})$, $g(\mathcal{G})\in G_{*}(\mathcal{G})$.
Let $e = \max\{e_1,e_2\}$, where $e_1, e_2$ are the extents of $g(\mathcal{F})$ and $g(\mathcal{G})$, respectively. It follows that $\mathcal{F}$ and $\mathcal{G}$ are $[e+1, n]$-normal. Under the assumption that $e > m$, we will derive a contradiction. Let $u$ be a positive integer with $t\leq u\leq e$ such that $g_{u}^{*}(\mathcal{F})\neq \emptyset$. Then we have the following cases.

If $u = t$, we claim $e = t$. Assume $e > t$ and choose $E \in g_{t}^{*}(\mathcal{F})$. The condition $e > t$ implies the existence of an integer $x \in [e] \setminus E$. There exists $E_1 \in g(\mathcal{F})$ with $E_1 \subseteq s_{x,e}(E)$, then $|E_1 \cap E| < t$, a contradiction. Thus, $e = t$. But then $t = e > m$ yields a contradiction. If $u = e$, then $g(\mathcal{F}) = \{[e]\}$ and $\mathcal{F} = \left\{ F \in \binom{[n]}{k_1} : [e] \subseteq F \right\}$, which contradicts $T \subseteq \mathcal{F}$ since $e > m$. Thus, $t < u < e$, and by Lemma~\ref{t-intersecting-cross-int-lemma}, we consider the following two cases.

\textbf{Case 1:} All four families $g_{u}^{*}(\mathcal{F})$, $g_{e+t-u}^{*}(\mathcal{F})$, $g_{e+1-u}^{*}(\mathcal{G})$ and $g_{u+1-t}^{*}(\mathcal{G})$ are non-empty. We then distinguish the following subcases based on the value of $u$.

\textbf{Case 1.1:} $u\neq e+t-u$.

Take $\mathcal{F}_1$ and $\mathcal{G}_1$ in   Table~\ref{tab:constructions}, 
and let $|g^{*}_u(\mathcal{F})| = |g^{*}_{e+1-u}(\mathcal{G})| = \alpha$ and $|g^{*}_{e+t-u}(\mathcal{F})| = |g_{u+1-t}^{*}(\mathcal{G})| = \beta$. Then the maximality of $|\mathcal{F}|+|\mathcal{G}|$ implies that
\begin{equation}\label{equ-1}
    \alpha\binom{n-e}{k_2+u-e} + \beta\binom{n-e}{k_1+(u+1-t)-e} \leq \alpha\binom{n-e}{k_1-u} + \beta\binom{n-e}{k_2-(u+1-t)}.
\end{equation}

Similarly, 
take $\mathcal{F}_2$ and $\mathcal{G}_2$ in Table~\ref{tab:constructions}. Then, the maximality of $|\mathcal{F}|+|\mathcal{G}|$ implies that
\begin{equation}\label{equ-2}
    \alpha\binom{n-e}{k_1-u+1} + \beta\binom{n-e}{k_2-(u-t)} \leq \beta\binom{n-e}{k_1+(u-t)-e} + \alpha\binom{n-e}{k_2+u-1-e}.
\end{equation}

Combining inequalities (\ref{equ-1}) and (\ref{equ-2}) and comparing the binomial coefficients—the detailed analysis is given in Appendix (Lemma~\ref{fact-1})—we obtain a contradiction.

\textbf{Case 1.2:} $u=e+t-u$. 

Take $\mathcal{F}_3$ and $\mathcal{G}_3$ in Table~\ref{tab:constructions}. 
Then the maximality of $|\mathcal{F}| + |\mathcal{G}|$ implies that
\begin{equation*}
    \binom{n-e}{k_1-u} \geq \binom{n-e}{k_2 - (e - u)}.
\end{equation*}
It follows that $n-e-(k_2-e+u)\geq k_1-u\geq k_2-e+u$, i.e., $k_1-k_2\geq 2u-e=t$, which is a contradiction.

\textbf{Case 2:} Only the families $g_{u}^{*}(\mathcal{F})$ and $g_{e+1-u}^{*}(\mathcal{G})$ are non-empty.

Take $\mathcal{F}_3$ and $\mathcal{G}_3$,  $\mathcal{F}_4$ and $\mathcal{G}_4$ in Table~\ref{tab:constructions},
then from the maximality of $|\mathcal{F}| + |\mathcal{G}|$, we obtain the inequalities:
\[
\binom{n-e}{k_1-u} \geq \binom{n-e}{k_2-(e-u)},
\quad
\binom{n-e}{k_2-(e+1-u)} \geq \binom{n-e}{k_1-(u-1)}.
\]

However, when $n > k_1 + k_2$, the following strict inequality holds:
\[
\binom{n-e}{k_1-u} \cdot \binom{n-e}{k_2-(e+1-u)} < \binom{n-e}{k_2-(e-u)} \cdot \binom{n-e}{k_1-(u-1)},
\]
which yields a contradiction.
\end{proof}

\subsection{The $[1,m]$-normal Property and Operations on Families $\mathcal{F}$ and $\mathcal{G}$}

In this subsection, we prove that under the conditions in our main theorem,  we can finally obtain two $[1,m]$-normal families after applying the operations described below. 

\begin{proof}[Proof of Lemma~\ref{[1,m]-normal}]
Assume that $\mathcal{F}, \mathcal{G}$ are $[\ell+1, m]$-normal but not $[\ell, m]$-normal. By Lemma~\ref{key-claim-frankl}(i), there exist a set $E_1 \subseteq [\ell]$ with $\ell \in E_1$ and $|E_1| = a$, and an interval $I(E_1) = [x_1, x_2]$ such that 
\[
E_1 \times \binom{[\ell+1,m]}{x} \subseteq u(\mathcal{F}) \quad \text{and} \quad E_1' \times \binom{[\ell+1,m]}{x+1} \cap u(\mathcal{F}) = \emptyset
\]
for all $x \in I(E_1)$. 
We repeat the following process for $\ell\le m-1$. By Lemma \ref{key-claim-frankl} (iv), we proceed by considering two cases based on the existence of $E_2$ to show that there is a contradiction or we can remove the critical set $E_1$ without changing the sum of the size of two families. When there are no critical sets $E_1$, the new families will be $[\ell,m]$-normal. Finally, $\mathcal{F}, \mathcal{G}$ can be modified to be $[1,m]$-normal by operations described in cases 1.2, 2.1.2, 2.2.1 and 2.2.2 below.

\subsection*{Case 1. There are no good pairs.}

In this case, we consider new cross-intersecting families: $\widetilde{\mathcal{F}}_1$, $\widetilde{\mathcal{G}}_1$, and $\widetilde{\mathcal{F}}_2$, $\widetilde{\mathcal{G}}_2$.
We recall the above specific constructions:
\begin{align*}
\widetilde{\mathcal{F}}_1&=\mathcal{F}\setminus\left\{\bigcup_{x \in I(E_1)} E_1\times\binom{[\ell+1,m]}{x}\times\binom{[m+1,n]}{k_1-a-x}\right\};\\
\widetilde{\mathcal{G}}_1&=\mathcal{G}\cup \left\{\bigcup_{x\in I(E_1)}F_1'\times\binom{[\ell+1,m]}{m-\ell-x}\times\binom{[m+1,n]}{k_2+a+x-m}\right\};
\end{align*}
and
\begin{align*}
\widetilde{\mathcal{F}}_2&=\mathcal{F}\cup \left\{\bigcup_{x \in I(E_1)} E_1'\times\binom{[\ell+1,m]}{x+1}\times\binom{[m+1,n]}{k_1-a-x}\right\}  \\
\widetilde{\mathcal{G}}_2&=\mathcal{G}\setminus \left\{\bigcup_{x\in I(E_1)}F_1\times\binom{[\ell+1,m]}{m-\ell-x-1}\times\binom{[m+1,n]}{k_2+a+x-m}\right\}.
\end{align*}

By the maximality of $|\mathcal{F}| + |\mathcal{G}|$, we obtain:
\begin{equation}\label{onlyE1-inq-1}
   \sum_{x \in I(E_1)} \binom{m-\ell}{x} \binom{n-m}{k_1-a-x} \geq \sum_{x \in I(E_1)} \binom{m-\ell}{m-\ell-x} \binom{n-m}{k_2+a+x-m},
\end{equation}
and
\begin{equation}\label{onlyE1-inq-2}
    \sum_{x \in I(E_1)} \binom{m-\ell}{x+1} \binom{n-m}{k_2+a+x-m} \geq \sum_{x \in I(E_1)} \binom{m-\ell}{x+1} \binom{n-m}{k_1-a-x}.
\end{equation}
Define $g(x) = \binom{n-m}{k_1-a-x} - \binom{n-m}{k_2+a+x-m}$ for $x \in I(E_1)$. We consider two subcases:

\textbf{Case 1.1.} $|I(E_1)| \geq 2$.

We claim that the sequence $(g(x) : x \in I(E_1))$, where $|I(E_1)| \geq 2$, satisfies the property: if $g(x) > 0$ then $g(x-1) > 0$, and if $g(x) < 0$ then $g(x+1) < 0$. The detailed verification of this property is provided in Appendix (Lemma~\ref{fact-2}). By Lemma~\ref{key-binomial-ineq}, inequalities~\eqref{onlyE1-inq-1} and~\eqref{onlyE1-inq-2} cannot both hold simultaneously, yielding a contradiction.

\textbf{Case 1.2.} $|I(E_1)| = 1$. 

Here $x_1 = x_2 = x$. Note that $g(x) = 0$ is equivalent to $k_1 - k_2 + m = 2(a + x)$. Now, only $E_1\times \binom{[\ell+1,m]}{x}$ violates the $[\ell,m]$-normality of $\mathcal{F}$, then 
we just delete $E_1\times \binom{[\ell+1,m]}{x}$ such that $E_1$ is not critical for $[\ell+1,m]$-normality of $\mathcal{F}$. Thus, we can assume that the families $\widetilde{\mathcal{F}}_1$ and $\widetilde{\mathcal{G}}_1$ are $[\ell, m]$-normal.

\subsection*{Case 2. There exists a good pair $(E_1,E_2)$.}

We further divide into subcases based on the value of $a$.

\textbf{Case 2.1.} $a \neq \ell + t - a$. 

In this case, consider $\widetilde{\mathcal{F}}_3$, $\widetilde{\mathcal{G}}_3$ and $\widetilde{\mathcal{F}}_4$, $\widetilde{\mathcal{G}}_4$. We recall the above specific constructions: 
\begin{align*}
  \widetilde{\mathcal{F}}_3&=\mathcal{F}\setminus \left\{\bigcup_{E_1\in G(E_1)}\bigcup_{x \in I(E_1)} E_1\times\binom{[\ell+1,m]}{x}\times\binom{[m+1,n]}{k_1-(a+x)}\right\} \\
  &\cup \left\{\bigcup_{E_2\in G(E_2)}\bigcup_{x\in I(E_1)}E_2'\times\binom{[\ell+1,m]}{m-\ell-x}\times\binom{[m+1,n]}{k_1+a+x-(t-1)-m}\right\},
\end{align*}
and
\begin{align*}   \widetilde{\mathcal{G}}_3&=\mathcal{G}\setminus\left\{\bigcup_{F_2\in G(E_2)}\bigcup_{x \in I(E_1)} F_2\times\binom{[\ell+1,m]}{x}\times\binom{[m+1,n]}{k_2-(a+x)+t-1}\right\} \\
&\cup\left\{\bigcup_{F_1\in G(F_1)}\bigcup_{x\in I(E_1)}F_1'\times\binom{[\ell+1,m]}{m-\ell-x}\times\binom{[m+1,n]}{k_2+a+x-m}\right\}.
\end{align*}
Similarly, we construct $\widetilde{\mathcal{F}}_4$ and $\widetilde{\mathcal{G}}_4$.  
Let $\alpha$ be the number of such $E_1$ (equivalently, $\alpha=|G(E_1)|$) and $\beta$ the number of such $E_2$ (equivalently, $\beta=|G(E_2)|$). By Lemma \ref{key-claim-frankl}, we deduce that $|G(F_1)|=\alpha$ and $|G(F_2)|=\beta$.

Then the maximality of $|\mathcal{F}|+|\mathcal{G}|$ implies:
\begin{align*}
&\sum_{x \in I(E_1)} \binom{m-\ell}{x} \left( \alpha \binom{n-m}{k_1-(a+x)} + \beta \binom{n-m}{k_2-(a+x)+t-1)} \right) \\
\geq&\sum_{x \in I(E_1)} \binom{m-\ell}{x} \left( \alpha \binom{n-m}{k_2+a+x-m} + \beta \binom{n-m}{k_1+a+x-(t-1)-m} \right),
\end{align*}
and
\begin{align*}
&\sum_{x \in I(E_1)} \binom{m-\ell}{x+1} \left( \alpha \binom{n-m}{k_2+a+x-m} + \beta \binom{n-m}{k_1+a+x-(t-1)-m} \right) \\
\geq&\sum_{x \in I(E_1)} \binom{m-\ell}{x+1} \left( \alpha \binom{n-m}{k_1-(a+x)} + \beta \binom{n-m}{k_2-(a+x)+t-1} \right).
\end{align*}
Define $f(x)=\alpha\left(\binom{n-m}{k_1-(a+x)}-\binom{n-m}{k_2+(a+x)-m}\right)+\beta\left(\binom{n-m}{k_2-(a+x)+t-1}-\binom{n-m}{k_1+(a+x)-(t-1)-m}\right)$ 
for $x \in I(E_1)$. We consider two subcases:

\textbf{Case 2.1.1.} $|I(E_1)| \geq 2$. 

We claim that the sequence $(f(x) : x \in I(E_1))$, where $|I(E_1)| \geq 2$, satisfies the following property: if $f(x) > 0$ for some $x \in I(E_1)$, then $f(x-1) > 0$; and if $f(x) < 0$ for some $x \in I(E_1)$, then $f(x+1) < 0$. The detailed verification of this property, involving explicit estimates, is provided in Appendix (Lemma~\ref{fact-3}). By Lemma~\ref{key-binomial-ineq}, the above two inequalities cannot hold simultaneously, yielding a contradiction.

\textbf{Case 2.1.2.} $|I(E_1)| = 1$. 

Let $x_1 = x_2 = x$. Then $f(x) = 0$. Take $\widetilde{\mathcal{F}}_1$ and $\widetilde{\mathcal{G}}_1$.  
Then the maximality of $|\mathcal{F}|+|\mathcal{G}|$ implies 
\begin{equation*}
    \binom{n-m}{k_1-(a+x)}\geq \binom{n-m}{k_2+(a+x)-m}.
\end{equation*}
If $\binom{n-m}{k_1-(a+x)} - \binom{n-m}{k_2+(a+x)-m} > 0$, then $f(x)=0$ means $\binom{n-m}{k_2-(a+1+x-t)} - \binom{n-m}{k_1+(a+1+x-t)-m} < 0$. It follows that $k_1 - k_2 < t - 1$, yielding a contradiction. Now the equality $f(x)=0$ holds if and only if $k_1 - k_2 = t - 1$ and $2(a + x) = m + t - 1$. Now, only $E_1\times \binom{[\ell+1,m]}{x}$ violates the $[\ell,m]$-normality of $\mathcal{F}$, then 
we just delete such terms such that $E_1$ is not critical for $[\ell+1,m]$-normality of $\mathcal{F}$. Thus, we can assume that the families $\widetilde{\mathcal{F}}_1$ and $\widetilde{\mathcal{G}}_1$ are $[\ell, m]$-normal.

\textbf{Case 2.2.} $a = \ell + t - a$. 

By Claim 2, $E \times \binom{[\ell+1,m]}{x_1}, E \times \binom{[\ell+1,m]}{x_2} \subseteq \mathcal{F}$ with $x_2 + x_1 = m - \ell - 1$. 
Take $\widetilde{\mathcal{F}}_5$ and $\widetilde{\mathcal{G}}_5$. We recall the above specific constructions: 
\begin{align*}
  \widetilde{\mathcal{F}}_5&=\mathcal{F}\setminus \left\{\bigcup_{E\in G(E)}\bigcup_{x \in I(E)} E\times\binom{[\ell+1,m]}{x}\times\binom{[m+1,n]}{k_1-(a+x)}\right\},\\  \widetilde{\mathcal{G}}_5&=\mathcal{G} \cup\left\{\bigcup_{F\in G(F)}\bigcup_{x\in I(E)}F'\times\binom{[\ell+1,m]}{m-\ell-x}\times\binom{[m+1,n]}{k_2+a+x-m}\right\}.
\end{align*}

Then the maximality of $|\mathcal{F}|+|\mathcal{G}|$ gives:
\begin{equation}\label{char-[1,m]-normal}
   \sum_{x \in I(E)} \binom{m-\ell}{x} \binom{n-m}{k_1-a-x} \geq \sum_{x \in I(E)} \binom{m-\ell}{m-\ell-x} \binom{n-m}{k_2+a+x-m}. 
\end{equation}

\textbf{Case 2.2.1.} There exists $x_0 \in I(E)$ such that $\binom{n-m}{k_1-a-x_0} = \binom{n-m}{k_2+a+x_0-m}$.

It follows that $2(a + x_0) = k_1 - k_2 + m$ since $n > k_1 + k_2$. Clearly, $2x_0 = k_1 - k_2 + m - (\ell + t) < m - \ell$ (i.e., $x_0 \leq \frac{m-\ell-1}{2}$). If $|I(E)| \geq 2$, let
$I(E) = \{x_0 - c, \ldots, x_0 - 1, x_0, x_0 + 1, \ldots, x_0 + d\}$ with $d \geq c$. For $i \leq c$, we have:
\begin{align*}
&\quad\binom{m-\ell}{x_0-i} \binom{n-m}{k_1-a-x_0+i} + \binom{m-\ell}{x_0+i} \binom{n-m}{k_1-a-x_0-i} \\
&< \binom{m-\ell}{m-\ell-x_0+i} \binom{n-m}{k_2+a+x_0-m-i} + \binom{m-\ell}{m-\ell-x_0-i} \binom{n-m}{k_2+a+x_0+i-m},
\end{align*}
since $\binom{m-\ell}{x_0+i} - \binom{m-\ell}{x_0-i} > 0$ and $\binom{n-m}{k_2+a+x_0+i-m} - \binom{n-m}{k_1-a-x_0-i} > 0$. Combined with $\binom{n-m}{k_1-a-x_0-j} < \binom{n-m}{k_2+a+x_0-m+j}$ for $j > 0$, this contradicts~\eqref{char-[1,m]-normal}. If $|I(E)| = 1$, then $x_0 = \frac{m-\ell-1}{2}$. Now $f(x) = 0$ holds if and only if $k_1 - k_2 = t - 1$ and $2(a + x) = m + t - 1$. Now, only $E\times \binom{[\ell+1,m]}{x}$ violates the $[\ell,m]$-normality of $\mathcal{F}$, then 
we just delete such terms such that $E$ is not critical for $[\ell+1,m]$-normality of $\mathcal{F}$. Thus, we can assume that the families $\widetilde{\mathcal{F}}_5$ and $\widetilde{\mathcal{G}}_5$ are $[\ell, m]$-normal. 

\textbf{Case 2.2.2.} $\binom{n-m}{k_1-a-x} \neq \binom{n-m}{k_2+a+x-m}$ for all $x \in I(E)$.

If $|I(E)| \geq 2$, let $x_0' = \lceil \frac{k_1 - k_2 + m}{2} \rceil - a < \frac{m-\ell}{2}$ (this inequality holds if and only if $k_1 - k_2 = t - 1$). Then $\binom{n-m}{k_1-a-x_0'} < \binom{n-m}{k_2+a+x_0'-m} = \binom{n-m}{k_1-a-x_0'+1}$, and hence $\binom{n-m}{k_1-a-x_0'-i} < \binom{n-m}{k_2+a+x_0'-m+i}$ for all $i > 0$. Let $I(E) = \{x_0' - c', \ldots, x_0' - 1, x_0', x_0' + 1, \ldots, x_0' + d'\}$ with $d' \geq c'$. For $i \leq c'$:
\begin{align*}
&\quad\binom{m-\ell}{x_0'-i} \binom{n-m}{k_1-a-x_0'+i} + \binom{m-\ell}{x_0'+i-1} \binom{n-m}{k_1-a-x_0'-(i-1)} \\
&< \binom{m-\ell}{m-\ell-x_0'+i} \binom{n-m}{k_2+a+x_0'-m-i} + \binom{m-\ell}{m-\ell-x_0'-(i-1)} \binom{n-m}{k_2+a+x_0'+i-1-m},
\end{align*}
since $\binom{m-\ell}{x_0'+i-1} > \binom{m-\ell}{x_0'-i}$ and $\binom{n-m}{k_2+a+x_0'+i-1-m} = \binom{n-m}{k_1-a-x_0'+i} > \binom{n-m}{k_1-a-x_0'-(i-1)} = \binom{n-m}{k_2+a+x_0'-m-i}$. Combined with $\binom{n-m}{k_1-a-x_0'-j} < \binom{n-m}{k_2+a+x_0'-m+j}$ for $j > 0$, this contradicts~\eqref{char-[1,m]-normal}.

If $|I(E)| = 1$, then $I(E) = \{\frac{m-\ell-1}{2}\}$ and $\binom{n-m}{k_1-a-\frac{m-\ell-1}{2}} \leq \binom{n-m}{k_2+a+\frac{m-\ell-1}{2}-m}$ since $k_1 - a - \frac{m-\ell-1}{2} \leq k_2 + a + \frac{m-\ell-1}{2} - m < n - m - k_1 + a + \frac{m-\ell-1}{2}$. The equality $f(x) = 0$ holds if and only if $k_1 - k_2 = t - 1$ and $2(a + x) = m + t - 1$. Similarly, we can assume that $\widetilde{\mathcal{F}}_5$ and $\widetilde{\mathcal{G}}_5$ are $[\ell, m]$-normal.

In summary, we show the normalization process as follows.

(i) The constructions $\widetilde{\mathcal{F}}_1,\widetilde{\mathcal{G}}_1$ yield $[\ell, m]$-normality under condition (1): $2(a + x) = k_1 - k_2 + m$; 

(ii) The constructions $\widetilde{\mathcal{F}}_3,\widetilde{\mathcal{G}}_3$ and $\widetilde{\mathcal{F}}_5,\widetilde{\mathcal{G}}_5$ yield $[\ell, m]$-normality under conditions (2) and (3): $2(a + x) = k_1 - k_2 + m$ and $k_1 - k_2 = t - 1$. 

Iterating this process, we eventually obtain $[1, m]$-normality. If all steps belong to Case 1, then $2(a + x) = k_1 - k_2 + m$ contradicts $2(a + x) \geq m + t$. If at least one step falls under Case 2, then $2(a + x) = k_1 - k_2 + m$ and $k_1 - k_2 = t - 1$. Thus, $\mathcal{F}$ and $\mathcal{G}$ are $[1, m]$-normal via suitable operations on sets $E \times \binom{[\ell+1,m]}{x} \times \binom{[m+1,n]}{k_1-(a+x)}$, which occur precisely when $2(a + x) = k_1 - k_2 + m$ and $k_1 - k_2 = t - 1$.
\end{proof}

\section{Proof of Theorem \ref{mainthm}}\label{proofofthm}

In this section, we prove our main result. By the properties of the shifting method, it is enough to prove the following lemma.

\begin{lemma}
Let $m \geq r \geq t$, $k_1 > k_2 > k_1 - t$, and $n \geq 2k_2 + t > k_1 + k_2 \geq r + k_2 \geq m + 1$, and $2r \geq m + t$ be integers.  Suppose $\mathcal{F}\subseteq\binom{[n]}{k_1}$ and $\mathcal{G}\subseteq\binom{[n]}{k_2}$ are left-compressed, cross-intersecting families, with $\mathcal{F}$ being $t$-intersecting and satisfying $\left\{F\in\binom{[n]}{k_1}:|F\cap [m]|\geq r\right\}\subseteq\mathcal{F}$, such that $|\mathcal{F}|+|\mathcal{G}|$ is maximized. Then the extremal families are:
\begin{enumerate}[(a)]
    \item $\mathcal{F} = \left\{F \in \binom{[n]}{k_1} : |F \cap [m]| \geq r\right\}$ and $\mathcal{G} = \left\{G \in \binom{[n]}{k_2} : |G \cap [m]| \geq m + 1 - r\right\}$.
    \item If $2(r - 1) = t - 1 + m$ and $k_1 - k_2 = t - 1$, then besides the extremal families in (a), we further have the following extremal families
    \[
    \mathcal{F} = \left\{F \in \binom{[n]}{k_1} : |F \cap [m]| \geq r\right\} \cup \mathcal{H}
    \]
    and
    \[
    \mathcal{G} = \left\{G \in \binom{[n]}{k_2} : |G \cap [m]| \geq m + 2 - r\right\} \cup \mathcal{K},
    \]
    where $\mathcal{H}=\mathcal{P}\times\binom{[m+1,n]}{k_1-r+1}$ satisfies that $\mathcal{P}\subseteq\binom{[m]}{r-1}$ 
    is $t$-intersecting, and 
    \[
    \mathcal{K}=\left(\binom{[m]}{m+1-r}\setminus\left\{[m]\setminus P:P\in \mathcal{P}\right\}\right) \times \binom{[m+1,n]}{k_2-(m+1-r)}.
    \]
\end{enumerate}
\end{lemma}

\begin{proof}
If $\mathcal{F}, \mathcal{G}$ are $[1, m]$-normal, then 
\[
\mathcal{F} = \left\{F \in \binom{[n]}{k_1} : |F \cap [m]| \geq a\right\}, \quad 
\mathcal{G} = \left\{G \in \binom{[n]}{k_2} : |G \cap [m]| \geq b\right\}
\]
with $2r \geq 2a \geq m + t$ and $a + b = m + 1$ by maximality of $|\mathcal{F}|+|\mathcal{G}|$. Increasing $a$ to $a + 1$ and decreasing $b$ to $b - 1$ yields a contradiction, since $n \geq 2k_2 + t > 2k_1 - t + (m + t - 2r) \geq 2(k_1 - m + b - 1) + m$ and $2b \leq 2(m + 1 - r)$ imply 
\[
\binom{m}{b-1} \binom{n-m}{k_2-(b-1)} > \binom{m}{a} \binom{n-m}{k_1-a} = \binom{m}{b-1} \binom{n-m}{k_1 - m + b - 1}.
\]
Hence, the maximality of $|\mathcal{F}|+|\mathcal{G}|$ implies that $a=r$.

If $\mathcal{F}, \mathcal{G}$ become $[1, m]$-normal after operations on sets $X \times \binom{[m+1,n]}{k_1 - a-x}$ with $|X| = a + x$, then $a + x = r - 1$ since $T \subseteq \mathcal{F}$. By Lemma \ref{[1,m]-normal}, we know that  $2(a + x) = k_1 - k_2 + m$ and $k_1 - k_2 = t - 1$. Denote these sets as $\mathcal{H}$. Notice that $\mathcal{H}$ has a disjoint decomposition as follows: 
\begin{align*}
  \mathcal{H}=\mathop{\dot{\bigcup}}\limits_{\substack{\ell\in A\subseteq [\ell],\\
  0\leq\ell\leq m-1}} A\times\binom{[\ell+1,m]}{x(A)}\times\binom{n-m}{k_1-|A|-x(A)}.
\end{align*}
where $x(A)$ is an integer and $\left(A'\times\binom{[\ell+1,m]}{x(A)+1}\right)\cap (\mathcal{H}\cap [m])=\emptyset$.

Then the maximality of $|\mathcal{F}|+|\mathcal{G}|$ implies
\[
\mathcal{F} = \left\{F \in \binom{[n]}{k_1} : |F \cap [m]| \geq r\right\} \cup \mathcal{H}
\]
and
\[
\mathcal{G} = \left\{G \in \binom{[n]}{k_2} : |G \cap [m]| \geq m + 2 - r\right\} \cup \mathcal{K},
\]
where $\mathcal{H}=\mathcal{P}\times\binom{[m+1,n]}{k_1-r+1}$ satisfies that $\mathcal{P}\subseteq\binom{[m]}{r-1}$ 
    is $t$-intersecting, and 
\[
\mathcal{K}=\left(\binom{[m]}{m+1-r}\setminus\left\{[m]\setminus P:P\in \mathcal{P}\right\}\right) \times \binom{[m+1,n]}{k_2-(m+1-r)}.
\]
\end{proof}

This completes the proof of Theorem~\ref{mainthm}.

\section{Appendix}

\begin{lemma}\label{fact-1}
If $n \geq 2k_2 + t > k_1 + k_2$ and $k_1 > k_2 > k_1 - t$, then the inequalities \eqref{equ-1} and \eqref{equ-2} cannot hold simultaneously.
\end{lemma}

\begin{proof}
Assume for contradiction that the \eqref{equ-1} and \eqref{equ-2} hold simultaneously. Multiplying inequalities \eqref{equ-1} and \eqref{equ-2}, the left-hand side becomes:
\begin{align*}
LHS = &\alpha^2 \binom{n-e}{k_1-u+1} \binom{n-e}{k_2+u-e} 
+ \beta^2 \binom{n-e}{k_1+(u+1-t)-e} \binom{n-e}{k_2-(u-t)} \\
+& \alpha\beta \left[ \binom{n-e}{k_1+(u+1-t)-e} \binom{n-e}{k_1-u+1} 
+ \binom{n-e}{k_2+u-e} \binom{n-e}{k_2-(u-t)} \right],
\end{align*}
while the right-hand side is:
\begin{align*}
RHS = &\alpha^2 \binom{n-e}{k_1-u} \binom{n-e}{k_2+u-1-e} 
+ \beta^2 \binom{n-e}{k_1+(u-t)-e} \binom{n-e}{k_2-(u+1-t)} \\
+& \alpha\beta \left[ \binom{n-e}{k_1+(u-t)-e} \binom{n-e}{k_1-u} 
+ \binom{n-e}{k_2+u-1-e} \binom{n-e}{k_2-(u+1-t)} \right].
\end{align*}

We now show that each term in the LHS strictly dominates its counterpart in the RHS. First, observe:
\begin{align*}
\binom{n-e}{k_1-u+1} \binom{n-e}{k_2+u-e} &>\binom{n-e}{k_1-u} \binom{n-e}{k_2+u-1-e}, \\
\binom{n-e}{k_1+(u+1-t)-e} \binom{n-e}{k_2-(u-t)} &> \binom{n-e}{k_1+(u-t)-e} \binom{n-e}{k_2-(u+1-t)}, \\
\binom{n-e}{k_1+(u+1-t)-e} \binom{n-e}{k_1-u+1} &> \binom{n-e}{k_1+(u-t)-e} \binom{n-e}{k_1-u}, \\
\binom{n-e}{k_2+u-e} \binom{n-e}{k_2-(u-t)} &> \binom{n-e}{k_2+u-1-e} \binom{n-e}{k_2-(u+1-t)}.
\end{align*}

These inequalities follow from the given conditions:
\begin{align*}
&n - e \geq k_1 + k_2+1 - e> (k_1 - u) + (k_2 + u - 1 - e), \\
&n - e \geq k_1 + k_2 +1- e > [k_1 + (u-t) - e] + [k_2 - (u+1-t)], \\
&n - e \geq 2k_1 - t+2 - e> [k_1 + (u-t) - e] + (k_1 - u), \\
&n - e \geq 2k_2 + t - e>(k_2+u-1-e)+[k_2-(u+1-t)].
\end{align*}

Now we conclude that $LHS > RHS$, contradicting the assumption that both \eqref{equ-1} and \eqref{equ-2} hold. This completes the proof.
\end{proof}

\begin{lemma}\label{fact-2}
If $n \geq 2k_2 + t > k_1 + k_2$ and $k_1 > k_2 > k_1 - t$, then the function 
\[
g(x) = \binom{n-m}{k_1-a-x} - \binom{n-m}{k_2+a+x-m}
\]
satisfies the following monotonicity property: if $g(x) > 0$ for some $x$, then $g(x-1) > 0$; and if $g(x) < 0$ for some $x$, then $g(x+1) < 0$.
\end{lemma}

\begin{proof}
We prove both conclusions by contradiction.

Suppose $g(x) > 0$ but $g(x-1) \leq 0$ for some $x$. Then:
\begin{align*}
\binom{n-m}{k_1-a-x} &> \binom{n-m}{k_2+a+x-m}, \\
\binom{n-m}{k_1-a-(x-1)}&\leq\binom{n-m}{k_2+a+(x-1)-m}.
\end{align*}

Multiplying these two inequalities yields:
\[
\binom{n-m}{k_1-a-x} \cdot \binom{n-m}{k_2+a+(x-1)-m} > \binom{n-m}{k_2+a+x-m} \cdot \binom{n-m}{k_1-a-(x-1)}.
\]

However, under the condition $n > k_1 + k_2$, we have the known binomial inequality:
\[
\binom{n-m}{k_2+a+x-m} \cdot \binom{n-m}{k_1-a-(x-1)} > \binom{n-m}{k_2+a+(x-1)-m} \cdot \binom{n-m}{k_1-a-x},
\]
which contradicts the previous inequality.

Suppose $g(x) < 0$ but $g(x+1) \geq 0$ for some $x$. Then:
\begin{align*}
\binom{n-m}{k_1-a-x} &< \binom{n-m}{k_2+a+x-m}, \\
\binom{n-m}{k_1-a-(x+1)} &\geq \binom{n-m}{k_2+a+(x+1)-m}.
\end{align*}

Multiplying these inequalities gives:
\[
\binom{n-m}{k_1-a-x} \cdot \binom{n-m}{k_2+a+(x+1)-m} < \binom{n-m}{k_2+a+x-m} \cdot \binom{n-m}{k_1-a-(x+1)}.
\]

But the known binomial inequality states:
\[
\binom{n-m}{k_1-a-x} \cdot \binom{n-m}{k_2+a+(x+1)-m} > \binom{n-m}{k_1-a-(x+1)}\cdot \binom{n-m}{k_2+a+x-m},
\]
which again leads to a contradiction.

Therefore, both conclusions hold as stated.
\end{proof}

\begin{lemma}\label{fact-3}
If $n \geq 2k_2 + t > k_1 + k_2$ and $k_1 > k_2 > k_1 - t$, then the function
$f(x) = \alpha\left(\binom{n-m}{k_1-(a+x)} - \binom{n-m}{k_2+(a+x)-m}\right) + \beta\left(\binom{n-m}{k_2-(a+1+x-t)} - \binom{n-m}{k_1+(a+1+x-t)-m}\right)$
satisfies: if $f(x) > 0$ then $f(x-1) > 0$, and if $f(x) < 0$ then $f(x+1) < 0$.
\end{lemma}

\begin{proof}
We prove the first implication by contradiction. Suppose $f(x) > 0$ but $f(x-1) \leq 0$. Then:
\begin{align*}
&\alpha\binom{n-m}{k_1-(a+x)} + \beta\binom{n-m}{k_2-(a+1+x-t)} > \alpha\binom{n-m}{k_2+(a+x)-m} + \beta\binom{n-m}{k_1+(a+1+x-t)-m}, \\
&\alpha\binom{n-m}{k_2+(a+x-1)-m} + \beta\binom{n-m}{k_1+(a+x-t)-m} \geq \alpha\binom{n-m}{k_1-(a+x-1)} + \beta\binom{n-m}{k_2-(a+x-t)}.
\end{align*}

Multiplying these inequalities, the left-hand side is: $LHS=$ 
\begin{align*}
&\quad\alpha^2 \binom{n-m}{k_1-(a+x)} \binom{n-m}{k_2+(a+x-1)-m} + \beta^2 \binom{n-m}{k_1+(a+x-t)-m} \binom{n-m}{k_2-(a+1+x-t)} \\
&+ \alpha\beta \left[ \binom{n-m}{k_1-(a+x)} \binom{n-m}{k_1+(a+x-t)-m} + \binom{n-m}{k_2+(a+x-1)-m} \binom{n-m}{k_2-(a+1+x-t)} \right],
\end{align*}
and the right-hand side is: $RHS =$
\begin{align*}
&\quad\alpha^2 \binom{n-m}{k_1-(a+x-1)} \binom{n-m}{k_2+(a+x)-m} + \beta^2 \binom{n-m}{k_1+(a+1+x-t)-m} \binom{n-m}{k_2-(a+x-t)} \\
&+ \alpha\beta \left[ \binom{n-m}{k_1-(a+x-1)} \binom{n-m}{k_1+(a+1+x-t)-m} + \binom{n-m}{k_2-(a+x-t)} \binom{n-m}{k_2+(a+x)-m} \right].
\end{align*}

However, the following binomial inequalities hold:
\begin{align*}
&\binom{n-m}{k_1-(a+x)} \binom{n-m}{k_2+(a+x-1)-m} < \binom{n-m}{k_1-(a+x-1)} \binom{n-m}{k_2+(a+x)-m}, \\
&\binom{n-m}{k_1+(a+x-t)-m} \binom{n-m}{k_2-(a+1+x-t)} < \binom{n-m}{k_1+(a+1+x-t)-m} \binom{n-m}{k_2-(a+x-t)}, \\
&\binom{n-m}{k_1-(a+x)} \binom{n-m}{k_1+(a+x-t)-m} < \binom{n-m}{k_1-(a+x-1)} \binom{n-m}{k_1+(a+1+x-t)-m}, \\
&\binom{n-m}{k_2+(a+x-1)-m} \binom{n-m}{k_2-(a+1+x-t)} < \binom{n-m}{k_2-(a+x-t)} \binom{n-m}{k_2+(a+x)-m}.
\end{align*}
These follow from the conditions:
\begin{align*}
&n - m \geq k_1 + k_2 + 1 - m > [k_1 - (a+x)]+ [k_2 + (a+x-1) - m], \\
&n - m \geq k_1+k_2+1-m>[k_1+(a+x)-t-m]+[k_2-(a+x)+t-1], \\
&n - m \geq 2k_1-t+2-m > [k_1-(a+x)]+[k_1+(a+x)-t-m], \\
&n - m \geq 2k_2 + t - m>[k_2+(a+x)-1-m]+[k_2-(a+x)+t-1].
\end{align*}

Now we have $LHS < RHS$, contradicting the assumption that $f(x) > 0$ and $f(x-1) \leq 0$.

The second implication (if $f(x) < 0$ then $f(x+1) < 0$) follows by a similar argument.
\end{proof}

\begin{lemma}\label{fact-4}
The inequalities (\ref{lem-ineq-3}) and (\ref{lem-ineq-4}) cannot hold simultaneously if $\bar{x}-x_1\geq 1$. And if $\bar{x}=x_1$, the two inequalities hold if and only if $2(a+\bar{x})=k_1-k_2+m$.
\end{lemma}
\begin{proof}
    Let $g(x) = \binom{n-m}{k_1-a-x} - \binom{n-m}{k_2+a+x-m}$, $h(x)=\binom{m-\ell}{x}g(x)$, $I=[x_1,\bar{x}]$, $c=k_2-k_1+t-1\geq 0$ and $\xi=\frac{k_1-k_2+m}{2}-a$. Then $g(x)>0$ $(h(x)>0)$ if $x<\xi$, and $g(x)<0$ $(h(x)<0)$ if $x>\xi$ as $n>k_1+k_2$. 

    Assume that the inequalities (\ref{lem-ineq-3}) and (\ref{lem-ineq-4}) hold simultaneously. If $x_1\geq c$, then the inequality (\ref{lem-ineq-4}) can be simplified as $\sum_{y\in I-c}\binom{m-\ell}{y+c+1}g(y)\leq 0$, and the inequality (\ref{lem-ineq-3}) is $\sum_{x\in I}\binom{m-\ell}{x}g(x)\geq 0$. Define $R(z)=\frac{\binom{m-\ell}{z+c+1}}{\binom{m-\ell}{z}}$ and $R(z)$ is strictly decreasing with respect to $z$. Then by the property of $g(x)$, $R(z)g(z)\geq R(\xi)g(z)$ for all $z$, and further
    \begin{equation*}
        \binom{m-\ell}{z+c+1}g(z)\geq R(\xi)\binom{m-\ell}{z}g(z).
    \end{equation*}
    It follows that $\sum_{y\in I-c}\binom{m-\ell}{y}g(y)\leq 0$. Notice that $\sum_{x\in I}\binom{m-\ell}{x}g(x)\geq 0$ implies that $\sum_{y\in I-c}\binom{m-\ell}{x}g(x)> 0$ if $c>0$, which is a contradiction. 
    
    If $x_1<c$, then the inequality (\ref{lem-ineq-4}) becomes $\sum_{z\in [x_1,c)}\binom{m-\ell}{z+1}g(z-c)+\sum_{z\in [c,\bar{x}]}\binom{m-\ell}{z+1}g(z-c)$. The first part $S_1$ is greater than $0$ as $g(x)>0$. The second part is $S_2:=\sum_{z\in [c,\bar{x}]}\binom{m-\ell}{z+1}g(z-c)=\sum_{y\in [0,\bar{x}-c]}\binom{m-\ell}{y+c+1}g(y)\geq R(\xi)\sum_{y\in [0,\bar{x}-c]}\binom{m-\ell}{y}g(y)$. We claim that $S_2\geq 0$. If $\bar{x}-c\leq \xi$, then the property of $h(x)\geq 0$ if $x\leq \xi$ implies that $S_2\geq 0$. If $\bar{x}-c> \xi$, the condition $\sum_{x\in I}h(x)\geq 0$ implies that $S_2\geq 0$. Now we complete the proof of the claim. Therefore, $S_1+S_2>0$ again leads a contradiction.

   The condition $c=0$ means that $k_1-k_2=t-1$. Then two inequalities (\ref{lem-ineq-3}) and (\ref{lem-ineq-4}) become $\sum_{x\in I}\binom{m-\ell}{x}g(x)\geq 0$ and $\sum_{x\in I}\binom{m-\ell}{x+1}g(x)\leq 0$. By Lemma~\ref{key-binomial-ineq} and Lemma~\ref{fact-2}, we obtain a contradiction if $|I|\geq 2$. If $\bar{x}=x_1$, $\binom{n-m}{k_1-(a+\bar{x})}=\binom{n-m}{k_2+a+\bar{x}-m}$ holds if and only if $2(a+\bar{x})=k_1-k_2+m$.

\end{proof}


\begin{thebibliography}{10}

\bibitem{ahlswede1997complete}
Ahlswede, Rudolf; Khachatrian, Levon H.
\newblock The complete intersection theorem for systems of finite sets.
\newblock {\em  European J. Combin.} \textbf{18} ({1997}), no. 2, 125-136.


\bibitem{DF2}
Deza, Michel-Marie; Frankl, Peter.
\newblock Erd\H{o}s-Ko-Rado theorem--22 years later.
\newblock{SIAM J. Algebraic Discrete Methods} \textbf{4} (1983), no. 4, 419-431.

\bibitem{ellis2022intersection}
Ellis, David.
\newblock Intersection problems in extremal combinatorics: theorems, techniques and questions old and new. 
\newblock{\em London Math. Soc. Lecture Note Ser.}, 481
Cambridge University Press, Cambridge, 2022, 115-173 


\bibitem{erdos1961intersection}
Erd\H{o}s, Paul; Ko, Chao; Rado, Richard.
\newblock Intersection theorems for systems of finite sets.
\newblock {\em Quart. J. Math. Oxford Ser.} (2) \textbf{12} (1961), 313–320.


\bibitem{frankl2020erdHos}
Frankl, Peter.
\newblock Erd{\H{o}}s--Ko--Rado Theorem for a restricted universe.
\newblock{\em Electron. J. Combin.} \textbf{27} (2020), no. 2, Paper No. 2.18, 5 pp.

\bibitem{frankl2016new}
Frankl, Peter.
\newblock New inequalities for cross-intersecting families.
\newblock{\em Mosc. J. Comb. Number Theory}
\textbf{6} (2016), no. 4, 27–32.

\bibitem{frankl1978} Frankl, Peter.
\newblock The Erd\H{o}s-Ko-Rado theorem is true for $n=ckt$. 
\newblock{\em Combinatorics (Proc. Fifth Hungarian Colloq., Keszthely, 1976)}, Voll. I, pp. 365-375. 
\newblock{\em Colloq. Math. Soc. J\'{a}nos Bolyai}, 18 North-Holland Publishing Co. Amsterdam-New York, 1978.


\bibitem{frankl2016invitation}
Frankl, Peter; Tokushige, Norihide.
\newblock Invitation to intersection problems for finite sets.
\newblock{\em J. Combin. Theory Ser. A} \textbf{144} (2016), 157–211.


\bibitem{hilton1977intersection}
Hilton, Anthony J. W.
\newblock An intersection theorem for a collection of families of subsets of a finite set.
\newblock {\em J. London Math. Soc.} (2) \textbf{15} (1977), no. 3, 369–376.

\bibitem{hilton1967some}
Hilton, Anthony J. W; Milner, Eric.
\newblock Some intersection theorems for systems of finite sets.
\newblock {\em Quart. J. Math. Oxford Ser.}(2) \textbf{18} (1967), 369–384.


\bibitem{k}
Katona, Gyula O.H.
\newblock A simple proof of the Erd\H{o}s-Chao Ko-Rado theorem. 
\newblock{\em J. Combinatorial Theory Ser. B} \textbf{13} (1972), 183-184.

\bibitem{k68}
Katona, Gyula O.H.
\newblock A theorem of finite sets.
\newblock{\em Theory of Graphs (Proc. Colloq., Tihany, 1966)}, pp. 187-207. 
Academic Press, New York-London, 1968. 

\bibitem{katonasurvey}
Katona, Gyula O.H.
\newblock Results on intersecting families of subsets, a survey. \newblock{\em New Trends in algebras and combinatorics}, 
World Scientific Publishing Co. Pte. Ltd., Hackensack, NJ, 2020, pp. 134-148.



\bibitem{k63} 
Kruskal, Joseph B.
\newblock The number of simplices in a complex.
\newblock{\em  Mathematical optimization techniques,} University of California Press, Berkeley-Los Angeles, Calif., 1963, pp. 251-278


\bibitem{li2025non}
Li, Anshui; Zhang, Huajun.
\newblock On non-empty cross-$t$-intersecting families.
\newblock{\em J. Combin. Theory Ser. A} \textbf{210} (2025), Paper No. 105960, 14 pp.

\bibitem{li2013intersecting}
Li, Wei-Tian; Chen, Bor-Liang; Huang, Kuo-Ching; Lih, Ko-Wei.
\newblock Intersecting $k$-uniform families containing all the $k$-subsets of a given set.
\newblock{\em Electron. J. Combin.} \textbf{20} (2013), no. 3, Paper 38, 11 pp.

\bibitem{wang2017intersecting}
Wang, Jun; Zhang, Huajun.
\newblock Intersecting $k$-uniform families containing a given family.
\newblock{\em Discrete Math.} \textbf{340} (2017), no. 2, 140–144.


\bibitem{W84}
Wilson, Richard M.
\newblock The exact bound in the Erd\H{o}s-Ko-Rado theorem.
\newblock{\em Combinatorica} \textbf{4} (1984), no. 2-3, 247–257.


\bibitem{wu2023refined}
Wu, Biao.
\newblock A refined result on cross-intersecting families. 
\newblock{\em Discrete Appl. Math.}
\textbf{339} (2023), 149–153.


\bibitem{wu2025result}
Wu, Yongjiang; Feng, Lihua; Li, Yongtao.
\newblock A result for hemi-bundled cross-intersecting families.
\newblock{\em Adv. in Appl. Math.}
\textbf{169} (2025), Paper No. 102912, 21 pp.


\bibitem{wu2025proof}
Wu, Yongjiang; Feng, Lihua; Li, Yongtao.
\newblock Proof of Frankl's conjecture on cross-intersecting families.
\newblock{\em J. Combin. Theory Ser. A}
\textbf{216} (2025), Paper No. 106062, 13 pp.



\bibitem{zhang2025conjecture}
Zhang, Huajun; Wu, Biao.
\newblock On a conjecture of Tokushige for cross-$t$-intersecting families.
\newblock{\em J. Combin. Theory Ser. B}
\textbf{171} (2025), 49–70.



\end{thebibliography}
\end{document}